\documentclass[11pt]{article}

\usepackage{latexsym,amsfonts,amsmath,amssymb}
\usepackage{extarrows}
\usepackage[english]{babel}
\usepackage[latin1]{inputenc}
\usepackage[T1]{fontenc}
\usepackage{graphics}
\usepackage[active]{srcltx}

\begin{document}

\newtheorem{duge}{Lemma}[section]
\newtheorem{rem}[duge]{Remark}
\newtheorem{prop}[duge]{Proposition}
\newtheorem{defi}[duge]{Definition}
\newtheorem{theo}[duge]{Theorem}
\newtheorem{nota}[duge]{Notation}
\newtheorem{cor}[duge]{Corollary}
\newtheorem{rappel}[duge]{Rappel}
\newtheorem{hypo}[duge]{Hypothesis}
\newcommand{\mcp}{\mathbb{P}}
\newcommand{\mce}{\mathbb{E}}
\newcommand{\taux}{\lambda_{b,k}}
\newcommand{\pn}{\mathcal{P}_{n}}
\newcommand{\mcn}{\mathbb{N}}
\newcommand{\mcr}{\mathbb{R}}
\newcommand{\pinf}{\mathcal{P}_{\infty}}
\newcommand{\partn}{\{ 1,\dots,n \}}
\newcommand{\sr}{\mathcal{C}}
\newcommand{\ur}{\mathcal{U}}
\newcommand{\srb}{\mathcal{S}^{\downarrow}}
\newcommand{\un}{\mathbf{1}}
\newcommand{\munf}{\overrightarrow{\mu_n}}
\newcommand{\shiftn}{\overset{\rightarrow n}{\Gamma}}
\newcommand{\freq}[1]{\ensuremath{|{#1}|^{\downarrow}}}
\newcommand{\pc}{\widehat{p}_{t,s}}

\title{Fragmentation of   compositions and intervals}

\author{
Anne-Laure Basdevant}
\date{}
\maketitle

\begin{center}
\it{Laboratoire de Probabilités et Modèles Aléatoires,\\
 Universit\'e Pierre et Marie Curie,\\ 175 rue du Chevaleret,
75013 Paris, France.}
\end{center}

\vspace*{0.8cm}

\begin{abstract}

The fragmentation processes of exchangeable partitions have already
been studied by several authors. In this paper, we examine rather
fragmentation of exchangeable compositions, that means partitions of
$\mcn$ where the order of the blocks counts. We will prove that such
a fragmentation is bijectively associated to an interval
fragmentation. Using this correspondence, we  then calculate  the Hausdorff dimension of certain random closed
set that arise in interval fragmentations and we study Ruelle's
 interval fragmentation.
\end{abstract}


\bigskip
\noindent{\bf A.M.S. Classification. } 60 J 25, 60 G 09.

\bigskip
\noindent{\bf e-mail. } Anne-Laure.Basdevant@ens.fr

\vspace*{0.8cm}

\section{Introduction}
Random fragmentations describe an object which splits as time
passes. Two types of fragmentation have received a special attention
 : fragmentation of partitions of $\mcn$ and mass-fragmentation,
i.e. fragmentation on the space $\srb=\{s_1\ge s_2\ge
\ldots\ge 0, \sum_i s_i\le 1\}$.
 Berestycki \cite{Berestycki02} has
proved that for each homogeneous fragmentation process of
exchangeable partitions, we can canonically associate a mass
fragmentation. More precisely, according to the work of Kingman
\cite{Kingman82}, we know that if $\pi=(\pi_1,\pi_2,\ldots)$ is an
exchangeable random partition of $\mcn$ (i.e. the distribution of
$\pi$ is invariant under finite permutation of $\mcn$),  the
asymptotic frequency of  block $\pi_i$,
$f_i=\lim_{n\rightarrow\infty}\frac{Card(\pi_i\cap
\{1,\ldots,n\})}{n}$, exists a.s. We denote by $(\freq{\pi_i})_{i\in
\mcn}$ the sequence $(f_i)_{i\in\mcn}$ after a decreasing
rearrangement. If $(\Pi(t),t\ge 0)$ is a fragmentation of
exchangeable partitions, then $(\freq{\Pi_i(t)}_{i\in \mcn},t\ge 0)$
is a mass fragmentation. Conversely, a fragmentation of exchangeable
partitions can be constructed from a mass fragmentation via a
"paintbox process".

One of our goal in this paper is to develop an analog theory for
fragmentations of exchangeable compositions and interval
fragmentations. The notion of composition structure has been
introduced by Gnedin \cite{Gnedin97} ; roughly speaking, it can be
thought of a partition where the order of the block counts. Gnedin
proved a  theorem analogous to Kingman's Theorem in the case of
exchangeable compositions : for each probability measure $P$ that
describes the law of a random exchangeable composition, we can find
a probability measure on the open subset of [0,1], such that $P$ can
be recovered via a "paintbox process". This is why it seems very
natural to look for a correspondence between fragmentations of
compositions and interval fragmentations.

The first part of this paper develops the relation between
probability laws of exchangeable compositions and laws of random open
subsets, and its extension to infinite measures. Then we prove that
there exists indeed a one to one correspondence between
fragmentation of compositions and interval fragmentations. The next
part  gives some properties and characteristics of these processes
and briefly presents how this theory can be extended to
time-inhomogeneous fragmentations and self-similar fragmentations.

We then turn our attention to the estimation of the Hausdorff
dimension of random closed sets which arise in an interval
fragmentation. Finally,
as an application of this theory, we study in Section \ref{Ruelle} a
well known interval fragmentation introduced by Ruelle
\cite{Basdevant05, Bertoinlegall00, Bolthausensznitman98, Ruelle87}
and we give a description of its semi-group of transition.

\vspace*{0.5cm}

\section{Exchangeable compositions and open subsets of $]0,1[$}
\subsection{Probability measures}\label{casfini}

In this section, we define  exchangeable compositions following
Gnedin  \cite{Gnedin97}, and
 recall some useful properties.\\
For $n \in \mcn$, let $[n]$ be the set of integers $\{1,\ldots,n\}$
and write $[\infty]=\mcn$.

\begin{defi}For $n\in \mcn$, a composition of $[n]$  is an ordered sequence of disjoint, non
empty subsets
of $[n]$, $\gamma=(A_1,\ldots,A_k)$, with $\cup A_i=[n]$.\\
We denote by $\sr_{n}$ the set of  composition of $[n]$.
\end{defi}

Let  $\rho_n : \sr_{n}\rightarrow \sr_{n-1}$ be the restriction of a
composition of $[n]$ to a composition of $[n-1]$ and let $\sr$ be
the projective limit of $(\sr_{n},\rho_{n})$. We endow $\sr$ with
the product topology, then it is a compact set.

We say  that a sequence $(P_n)_{n\in \mcn}$ of measure on
$(\sr_n)_{n\in \mcn}$ is a consistent sequence of measures if, for
all $n\ge 2$, $P_{n-1}$ is the image of $P_n$ by the projection
$\rho_n$, i.e.,
 for  all  $\gamma\in \sr_{n-1}$, we have
$$P_{n-1}(\Gamma_{n-1}=\gamma)=\sum_{\gamma'\in \sr_n : \rho_n(\gamma')=\gamma}{P_n(\Gamma_n=\gamma)}.$$

By Kolmogorov theorem, such a sequence
 $(P_n)_{n\in \mcn}$
determines the law of a random composition of $\mcn$.

In the sequel, for $n\in \mcn\cup\{ \infty\}$,  $\gamma \in
\sr_n$ and
$A\subset [n]$, $\gamma_A$ will denote the restriction of $\gamma$ to $A$.
Hence, for  $m\le n$ and, $\gamma_{[m]}$ will denote the restriction of $\gamma$ to $[m]$.

A random composition $\Gamma$ of   $\mcn$ is called exchangeable if
for all $n \in \mcn$, for every permutation $\sigma$ of $[n]$ and
for all $\gamma\in \sr_{n}$, we have :
$$\mcp(\Gamma_{[n]}=\gamma)=\mcp(\sigma(\Gamma_{[n]})=\gamma),$$
where $\sigma(\Gamma_{[n]})$ the image of the composition $\Gamma_{[n]}$ by $\sigma$.
Hence, given an exchangeable random composition  $\Gamma$, we can
associate a function defined on finite sequences of $\mcn$ by
$$\forall k\in \mcn, \forall n_1,\ldots,n_k\in \mcn^k, p(n_1,\ldots,n_k)=\mcp(\Gamma_{[n]}=(B_1,\ldots,B_k)),$$
with $|B_i|=n_i$ and $n_1+\ldots+n_k=n$. This function determines
the law of $\Gamma$ and is called the exchangeable composition
probability function (ECPF) of $\Gamma$.

\begin{nota}
Let $\Gamma$ be a composition of $\mcn$. For $i,j\in \mcn^2$, we will use
the following notation :
\begin{itemize}
   \item $i\sim j$, if $i$ and $j$ are in the same block.
  \item $i\prec j$, if the block containing $i$ is before the  block containing $j$.
  \item  $i\succ j$, if the block containing $i$ is after the block containing $j$.
\end{itemize}
\end{nota}

\begin{defi} \label{paintbox}
Let $U$ be an open subset of $[0,1]$. We construct a random
composition of $\mcn$ in the following way :\\
Let us draw $(X_i)_{i\in \mcn}$ iid random variables with uniform
law on $[0,1]$. Then we use the following rules :
\begin{itemize}
   \item $i\sim j$, if $i=j$ or if $X_i$ and $X_j$ belong to the same component interval of  $U$.
  \item $i\prec j$,  if $X_i$ and $X_j$ do not belong to the same component interval of  $U$ and   $X_i<X_j$.
  \item  $i\succ j$,  if $X_i$ and $X_j$ do not belong to the same component interval of  $U$ and   $X_i>X_j$.
\end{itemize}
This defines a probability measure on $\sr$ that we shall denote
$P^{U}$ ;
 the marginal of $P^{U}$ on $\sr_n$ will be denoted by
 $P^{U}_{n}$.
 If $\nu$ is a
probability measure on $\ur$, we denote by $P^\nu$ the law on $\sr$
which marginals are :
$$P^{\nu}_{n}(\cdot)=\int_{\ur}{P^{U}_{n}(\cdot)d\nu(U)}.$$
\end{defi}

Let $\ur$ be the set of open subset of $]0,1[$. For $U\in \ur$,
let
$$\chi_U(x)=\min\{|x-y|,y\in U^{c}\}, \; x\in[0,1],$$ where
$U^{c}=[0,1]\backslash U$. We define also a distance on  $\ur$ by
:
$$d(U,V)=||\chi_U-\chi_V||_{\infty}.$$
It will be convenient to use the notation $\un=]0,1[$.  The
composition of $\sr_n$ (resp. $\sr$) with a single non empty block will be denoted by
$\un_{n}$ (resp. $\un_{\mcn}$) and we will write $\sr_n^{*}$ for $\sr_n\backslash
\{\un_{n}\}$.

Let us recall here two useful theorems from  Gnedin
\cite{Gnedin97} :
\begin{theo}~\cite{Gnedin97} \label{gnedinouv}\label{gnedin}
Let $\Gamma$ be an exchangeable random composition of $\mcn$, $\Gamma_{[n]}$
its restriction to $[n]$. Let $(n_1,\ldots,n_k)$ be the sequence of
the block sizes of $\Gamma_{[n]}$ and $n_0=0$. Define
 $U_n \in \ur$ by :
$$U_n=\bigcup_{i=1}^{k}\left]\frac{n_{i-1}}{n},\frac{n_i}{n}\right[.$$
Then $U_n$ converges almost surely to a random element $U\in \ur$.
The conditional law  of $\Gamma$ given $U$ is $P^U$.\\
As a consequence, if $P$ be an exchangeable probability  measure on $\sr$,
then there exists a unique probability
measure $\nu$ on $\ur$ such that $P=P^{\nu}$.
\end{theo}

Hence, for each exchangeable composition $\Gamma$, we can associate
an
element of $\ur$ which we will call  asymptotic open set of $\Gamma$ and  denote $U_{\Gamma}$.\\
We shall also write $|\Gamma|^{\downarrow}$  for the decreasing sequence
of the lengths of the interval components of $U_{\Gamma}$. More generally,
for $U\in \ur$, $|U|^{\downarrow}$ will be  the decreasing sequence
of the interval component lengths of $U$.

Let us notice that this theorem is the analogue of Kingman's Theorem
for the representation of exchangeable partitions. Actually, let
$\mathbf{Q}$ be an exchangeable probability measure on $\pinf$, the
set of partition of $\mcn$ and let $\pi$ be a partition with law
$\mathbf{Q}$.

Kingman \cite{Kingman82} has proved that each block of $\pi$  has
almost surely a frequency, i.e. if
 $\pi=(\pi_1,\pi_2,\ldots)$, then
 $$\forall i \in \mcn \hspace*{1cm}
f_i=\underset{n \rightarrow \infty}{\lim}\frac{\sharp\{\pi_i\cap
[n]\}}{n} \hspace*{1cm}\mbox{exists } \mathbf{Q} \mbox{-a.s.}$$ One
calls
 $f_i$  the frequency of the block $\pi_i$.
Therefore, for all exchangeable random partitions, we can associate
a probability on $\srb=\{s=(s_1,s_2,\ldots), s_1\ge s_2\ge \ldots\ge
0, \sum_i{s_i}\ge 1\}$ which will be the law of the decreasing
rearrangement of the sequence of the partition frequencies.

Conversely, given a law $\tilde{\nu}$ on $\srb$, we can construct an
exchangeable random partition whose law  of its  frequency  sequence
is $\tilde{\nu}$ (cf. \cite{Kingman82}) : we pick $s\in \srb$ with
law $\tilde{\nu}$ and we draw a sequence of independent random
variables $U_i$ with uniform law on $[0,1]$. Conditionally on $s$,
two integers $i$ and $j$ are in the same block of $\Pi$ iff there
exists an integer $k$ such that $\sum_{l=1}^{k}{s_l}\le U_i<
\sum_{l=1}^{k+1}{s_l} $ and $\sum_{l=1}^{k}{s_l}\le U_j<
\sum_{l=1}^{k+1}{s_l} $. We denote by $\rho_{\tilde{\nu}}$ the law
of this partition (and   by a slight abuse of notation, $\rho_{u}$
denotes the law of the partition obtained with
$\tilde{\nu}=\delta_u$). Kingman's representation Theorem states
that any  exchangeable random partition can be constructed in this
way.

Let $\wp_1$  be the  canonical application from the set of
composition $\sr$ to the set of partition $\pinf$ and $\wp_2$ the
application from the set $\ur$ to the set $\srb$ which associates to
an element $U$ of $\ur$ the decreasing sequence $\freq{U}$.
To sum up, we have  the following diagram between probability
measures on $\pinf$, $\sr$, $\srb$, $\ur$ :

$$
\begin{array}{ccc}
 (\sr,P^{\nu})&\xleftrightarrow{\mbox{ \tiny{Gnedin} }}& (\ur,\nu)\\
\mbox{\tiny{$\wp_1$}}\Big{\downarrow } & & \mbox{\tiny{$\wp_2$}}\Big{\downarrow }  \\
(\pinf,\rho_{\tilde{\nu}}) &
\xleftrightarrow{\mbox{\tiny{Kingman}}}& (\srb,\tilde{\nu}).
 \end{array}
$$

\subsection{Representation of infinite measures on
$\sr$}\label{mesure}

 In this section, we show how  Theorem \ref{gnedin} can be extended
to the case of an infinite measure $\mu$ on $\sr$ such that :
\begin{itemize}
\item $\mu$ is  exchangeable.
 \item $\mu(\un_{\mcn})=0.$
 \item For all $n \in \mcn$,  $\mu(\{\gamma \in \sr,
\gamma_{[n]}\neq \un_n \})<\infty$.
\end{itemize}

A measure on $\sr$ fulfilling this three properties will be called a
"fragmentation measure". We will  see in the sequel that such a
measure can always be associated to a fragmentation process and
conversely.

We will prove that we can decompose every fragmentation measure
$\mu$ in two measures, one  characterizing $\mu$ on the compositions
with asymptotic open set $U_{\gamma}=]0,1[$, and the other on the
complementary event. The measure on the event  $U_{\gamma}=]0,1[$ is
called erosion measure and the  measure on the event
$U_{\gamma}\neq]0,1[$ is called dislocation measure.

\begin{defi}
A measure $\nu$ on $\ur$ is called a dislocation measure if :
$$\nu(\un)=0, \hspace*{1.5cm} \int_{\ur}{(1-s_1)\nu(dU)}<\infty,$$
where $s_1$ is the length of the largest interval component of
$U$.
\end{defi}

In the sequel, for any $\nu$  measure on $\ur$, we define the
measure $P^{\nu}$  on $\sr$  by
$$P^{\nu}(d\gamma)=\int_{\ur}{P^{U}(d\gamma)d\nu(U)}.$$
Notice that if  $\nu$ is a dislocation measure, then
$P^{\nu}(d\gamma)$
 is a fragmentation measure.
In fact, the measure $P^{\nu}$ is exchangeable  since $P^{U}$ is
an exchangeable measure.\\
For $U\neq \un$, we have $P^{U}(\un_{\mcn})=0$, and as $\nu(\un)=0$,
we have also
 $P^{\nu}(\un_{\mcn})=0$.\\
We now have to check that  $P^{\nu}(\{\gamma \in \sr, \gamma_{[n]}\neq
\un_n \})<\infty$ for all $n \in \mcn$. Let us fix $U\in \ur$.
Set $\freq{U}=s=(s_1,s_2,\ldots)$.

$$P^{U}(\{\gamma \in \sr, \gamma_{[n]}\neq
\un_n \})=1-\sum_{i=1}^{\infty}{s_i^n}\le 1-s_1^n\le n(1-s_1)$$
and so $P^{\nu}(\{\gamma \in \sr, \gamma_{[n]}\neq \un_n
\})<\infty$.

\vspace*{0.5cm}

We can now state the following theorem :
\begin{theo}\label{disc}
Let $\epsilon_i$  be the composition of $\mcn$,
$(\{i\},\{\mcn\setminus\{i\}\})$ and
$\epsilon=\sum_{i}{\delta_{\epsilon_i}}.$ Let
$\epsilon'_i$  be the composition of $\mcn$,
$(\{\mcn\setminus\{i\}\},\{i\})$ and
$\epsilon'=\sum_{i}{\delta_{\epsilon'_i}}.$
These are  two exchangeable measures on $\sr$. \\
If $\mu$ is a fragmentation measure, there exists $c_l\ge 0$,
$c_r\ge 0$ and a dislocation measure  $\nu$ such that :
$$\mu=c_l\epsilon+c_r\epsilon'+P^{\nu}.$$
Besides, the restriction of $\mu$ to $\{\Gamma \in \sr ,\;
U_{\Gamma}=\mathbf{1}\}$ is $c_l\epsilon+c_r\epsilon'$ and the
restriction to  $\{\Gamma \in \sr ,\; U_{\Gamma}\neq\mathbf{1}\}$ is
$P^{\nu}$.
\end{theo}

Recall that in the case of fragmentation measure on partitions, Bertoin \cite{CoursBertoin03} proved the following result :\\
Let $\tilde{\epsilon_i}$  be the partition of $\mcn$,
$\big\{\{i\},\{\mcn\setminus\{i\}\}\big\}$ and   define the measure $\tilde{\epsilon}=\sum_{i}{\delta_{\tilde{\epsilon_i}}}$.
Let $\tilde{\mu}$ be an exchangeable measure on $\pinf$ such that $\mu(\un_{\mcn})=0$ and $\tilde{\mu}(\pi \in \pinf, \pi_n\neq \un_n)$ is finite for all
$n\in \mcn$. Then there exists  a measure  $\tilde{\nu}$ on $\srb$ such that $\tilde{\nu}(1)=0$ and $\int_{\srb}(1-s_1)\nu(ds)$, and a nonnegative number
$c$ such that :
$$\tilde{\mu}=\rho_{\tilde{\nu}}+c\tilde{\epsilon}.$$

Notice that Theorem \ref{disc} is
 an analogous decomposition as in the case of fragmentation measure on compositions,
except that, in this case, there is two coefficients of erosion, one
characterizing the left erosion and the other the right erosion.

\vspace*{0.3cm}

Proof. We adapt a proof    due to Bertoin
 ~\cite{CoursBertoin03} for the exchangeable partition to our case.\\
Set $n\in \mcn$. Set $\mu_n=\un_{\{\Gamma_{[n]}\neq \un_{n}\}}\mu$,
therefore $\mu_n$ is a finite measure. Let $\munf$ be the image of
$\mu_n$ by the $n$-shift, i.e. :
$$i\overset{\shiftn}{\prec} j\Leftrightarrow
i+n\overset{\Gamma}{\prec}j+n, \hspace*{0.7cm}
 i\overset{\shiftn}{\sim}
j\Leftrightarrow i+n\overset{\Gamma}{\sim}j+n, \hspace*{0.7cm}
i\overset{\shiftn}{\succ} j\Leftrightarrow
i+n\overset{\Gamma}{\succ}j+n.$$

Then  $\munf$ is exchangeable since $\mu$ is and furthermore, it
is finite measure. So, we can apply  Theorem \ref{gnedin} :
$$\exists\; !\; \nu_n \mbox{ finite measure on } \ur \mbox{ such that }
\munf(d\gamma)=\int_{\ur}{P^U(d\gamma)\nu_n(dU)}.$$

According to  theorem \ref{gnedinouv}, since $\munf$ is an
exchangeable finite measure,
 $\munf$-almost  every composition has
an asymptotic open set and so $\mu_n$-almost every  composition
has also an asymptotic open set, and as  $\mu=\lim\uparrow
\mu_n$, $\mu$-almost every composition has also an asymptotic open set.\\
Besides we have :
$$\mu_n(n+1\nsim n+2\;|\;U_\Gamma=U)=\munf(1\nsim 2\;|\;U_\Gamma=U)=P^U(1\nsim 2)=1-\sum{s_i²}\ge 1-s_1.$$
So $$\mu_n(n+1\nsim n+2)\ge \int (1-s_1)\nu_n(dU).$$ Set
$\nu=\lim_{n \rightarrow \infty} \uparrow \nu_n$. Since
$$\mu_n(n+1\nsim n+2)\le \mu(n+1\nsim n+2) \le  \mu(1\nsim
2)<\infty,$$ we have $$\int (1-s_1)\nu(dU)<\infty.$$ Hence $\nu$ is
a dislocation measure. Set $\gamma_k \in \sr_k$.
\begin{eqnarray*}
\mu(\Gamma_{[k]}=\gamma_k,U_{\Gamma}\neq \un)&=&\lim_{n\rightarrow \infty}\mu(\Gamma_{[k]}=\gamma_k,\Gamma_{\{k+1,\ldots,k+n\}}\neq \un_{n},
U_{\Gamma} \neq \un)\\
&=&\lim_{n\rightarrow \infty}\mu(\shiftn_{[k]}=\gamma_k,\Gamma_{[n]}\neq \un_{n}, U_{\Gamma} \neq \un)\\
&=&\lim_{n\rightarrow \infty}\munf(\Gamma_{[k]}=\gamma_k, U_{\Gamma} \neq \un)\\
&=&\int_{\sr^{*}}{P^U(\Gamma_{[k]}=\gamma_k)\nu(dU)}.
\end{eqnarray*}
Thus we have
$$\mu(\;\cdot\;,U_{\gamma}\neq \un)=\int P^U(\,\cdot\,)\nu(dU).$$

We now have to study $\mu$ on the event $\{ U_{\gamma}=\un\}$.\\
Let $\tilde{\mu}$ be $\mu$ restricted to $\{1\nsim 2,
U_{\gamma}=\un\}$. Let $\overset{\rightarrow}{\tilde{\mu}}$ be the
image of $\tilde{\mu}$ by the 2-shift.
 The measure $\overset{\rightarrow}{\tilde{\mu}}$ is  finite and exchangeable and its asymptotic open set is almost surely $\un$, so
 $\overset{\rightarrow}{\tilde{\mu}}=a \delta_{\un}$ where $a$ is a nonnegative number.\\
So
$\tilde{\mu}=c_1\delta_{\gamma_1}+\ldots+c_{10}\delta_{\gamma_{10}}$
where $\gamma_1,\ldots,\gamma_6$ are the six possible compositions build from the blocks $\{1\}$, $\{2\}$, $\mcn\backslash \{1,2\}$,\\
$\gamma_7=(\{1\}, \mcn \backslash \{1\})$, $\gamma_8=(\{2\}, \mcn \backslash \{2\})$,\\
$\gamma_9= (\mcn \backslash \{1\},\{1\})$, $\gamma_{10}= (\mcn
\backslash \{2\}, \{2\})$. We must have $c_1=\ldots=c_6=0$, for
otherwise, by exchangeability, we would have $\mu(\{1\},\{n\}, \mcn
\backslash \{1,n\})=c>0$ and this would yield $\mu(S_2^{*})=\infty$.
By exchangeability, we also have $c_7=c_8$ and $c_9=c_{10}$ and so,
by exchangeability,
$$\mu \un_{\{ U_{\gamma}=\un\}}=c_l\sum_i{\delta_{\epsilon_i}}+c_r\sum_i{\delta_{\epsilon'_i}}. \Box$$

As in section \ref{casfini}, we can now draw a diagram between
fragmentation measures on $\sr$ and $\pinf$ and dislocation measures
on $\ur$ and $\srb$. Let us recall that $\wp_1$ is the canonical
projection of $\sr$ to $\pinf$, and denote
$q:(\ur,\mcr_+,\mcr_+)\mapsto (\srb,\mcr_+)$ the application defined
by $q(U,a,b)=q(\freq{U},a+b)$. Then we have the following diagram :

$$
\begin{array}{ccc}
 (\sr,\mu)&\xleftrightarrow{\mbox{\tiny{ Theorem \ref{disc} }}}& \big(\ur,(\nu,c_l,c_r)\big)\\
\mbox{\tiny{$\wp_1$}}\Big{\downarrow } & &\mbox{\tiny{$q$}}\Big{\downarrow }  \\
(\pinf,\tilde{\mu})
&\xleftrightarrow{\quad\mbox{\tiny{Bertoin}\qquad}} &
\big(\srb,(\tilde{\nu},c_l+c_r)\big).
 \end{array}
$$

It remains to prove that $\tilde{\mu}=\rho_{\tilde \nu}+ (c_l+c_r)\tilde{\epsilon}$.
 Set $\tilde{\mu}=\rho_{\overline{\nu}}+c\tilde{\epsilon}$.
 Since $\tilde{\mu}$ is the image by $\wp_1$ of $\mu$, we have
 $$\tilde{\mu}(\tilde{\epsilon}_1)=\mu(\epsilon_1)+\mu(\epsilon'_1) \mbox{ and then } c=c_r+c_l.$$

 Let us fix $n\in \mcn$ and
 $\pi\in \pn\backslash\{\un_n\}$. Set $A=\{\gamma\in \sr_n, \wp_1(\gamma)=\pi\}$.
Remark now that for all $U,V \in \ur$ such that
$\freq{U}=\freq{V}$, we have $P^{U}(A)=P^V(A)$. Moreover we have
$P^{U}(A)=\rho_s(\pi)$ if $s=\freq{U}$. So
$$P^{\nu}(A)=\int_{\srb}P^U(A)\nu(U,\freq{U}=ds)=\int_{\srb}\rho_s(\pi)\tilde{\nu}(ds)=
\rho_{\tilde{\nu}}(\pi).$$
We get
$$\mu(A)=P^{\nu}(A)+c_l\epsilon(A)+c_r\epsilon'(A)=
\rho_{\tilde{\nu}}(\pi)+(c_l+c_r)\tilde{\epsilon}(A)
=\rho_{\overline{\nu}}(\pi)+(c_l+c_r)\tilde{\epsilon}(A)=\tilde{\mu}(\pi).$$
So we deduce that $\overline{\nu}=\tilde{\nu}$. $\Box$

\vspace*{0.3cm}

\section{Fragmentation of  compositions and interval fragmentation}
\subsection{Fragmentation of compositions}
\begin{defi}
Let us fix
 $n\in \mcn$ and $\gamma \in \sr_n$ with
$\gamma=(\gamma_1,\ldots,\gamma_k)$. Let
$\gamma^{(.)}=(\gamma^{(i)},i \in
\partn)$ with $\gamma^{(i)} \in \sr_n$ for all $i$. Set
$m_i=\min \gamma_i$. We denote $\tilde{\gamma}^{(i)}$ the
restriction of $\gamma^{(m_i)}$ to $\gamma_i$. So
$\tilde{\gamma}^{(i)}$ is a composition of $\gamma_i$. We consider
now
 $\tilde{\gamma}=(\tilde{\gamma}^{(1)},\ldots,\tilde{\gamma}^{(k)}) \in \sr_n$.\\
 We denote by $FRAG(\gamma,\gamma^{(.)})$ the composition $\tilde{\gamma}$.
 If $\gamma^{(.)}$ is a sequence of i.i.d. random variables with law
 $p$, $p$-FRAG$(\gamma,\cdot)$ will denote the law of FRAG$(\gamma,\gamma^{(.)}).$
\end{defi}

We remark then that the operator $FRAG$ has some useful property.
First, we have that $FRAG(\gamma,\mathbf{1}^{(.)})=\gamma$.
Furthermore, the fragmentation operator is compatible with the
restriction i.e. for every $n'\le n$ :
$$FRAG(\gamma,\gamma^{(.)})_{[n']}= FRAG(\gamma_{[n']},\gamma^{(.)}).$$
Besides, the operator $FRAG$ preserves the exchangeability. More
precisely, let   $(\gamma^{(i)},i\in \partn)$ be a sequence of
random compositions which is  \emph{doubly exchangeable}, i.e. for
each $i$, $\gamma^{(i)}$ is an exchangeable composition, and
moreover, the sequence
 $(\gamma^{(i)},i\in \partn)$ is also exchangeable. Let  $\gamma$ be an exchangeable composition of
$\sr_n$ independent of $\gamma^{(\cdot)}$. Then
$FRAG(\gamma,\gamma^{(.)})$ is an exchangeable composition. Let us
prove this property. Let us fix a permutation $\sigma$ of $[n]$. We
shall prove that
$$FRAG(\gamma,\gamma^{(.)})\overset{law}=\sigma(FRAG(\gamma,\gamma^{(.)})).$$
Let $k$ be the number of blocks of $\gamma$ and denote by
$m_1,\ldots,m_k$  the minimums of $\gamma_1,\ldots,\gamma_k$. Let
define now $m'_1, \ldots, m'_k$ the minimums of
$\sigma(\gamma_1),\ldots, \sigma(\gamma_k)$. Define now
$\gamma'^{(\cdot)}=(\gamma'^{(i)},i\in \partn)$ by
$$\gamma'^{(m'_i)}=\sigma(\gamma^{(m_i)}) \mbox{ for } 1\le i \le k$$
$$\gamma'^{(j)}=\sigma(\gamma^{(f(j))}) \mbox{ for } j \in \{1,\ldots,
n\}\setminus\{m'_i, 1\le i \le k\},$$ where $f$ is the increasing
bijection from $\{1,\ldots, n\}\setminus\{m'_i, 1\le i \le k\}$ to
$\{1,\ldots, n\}\setminus\{m_i, 1\le i \le k\}$. We get
$$\sigma(FRAG(\gamma,\gamma^{(.)}))=FRAG(\sigma(\gamma),\gamma'^{(.)}).$$
Since $\sigma(\gamma)\overset{law}= \gamma$ and
$\gamma'^{(.)}\overset{law}= \gamma^{(.)}$ and $\gamma'^{(.)}$
remains independent of $\gamma$, we get
$$FRAG(\sigma(\gamma),\gamma'^{(.)})\overset{law}{=}FRAG(\gamma,\gamma^{(.)}). \, \Box$$
We can now define the notion of exchangeable fragmentation process
of compositions.

\begin{defi}
Let us fix $n \in \mcn$ and let $(\Gamma_n(t),t\ge 0)$ be a Markov
process on $\sr_n$ which is continuous in
probability.\\
  We call $\Gamma_n$  an exchangeable fragmentation process of  compositions if :
 \begin{itemize}
 \item $\Gamma_n(0)=\un_{n}$ a.s.
\item  Its semi-group is described in the following way : there
exists a family of probability measures on exchangeable
compositions
 $(P_{t,t'},t\ge 0,t'> t)$ such that for all
$ t\ge 0,t^{'}> t$ the conditional law of $\Gamma_n(t^{'})$ given
 $\Gamma_n(t)=\gamma$ is the law of $P_{t,t'}$-FRAG$(\gamma,\gamma^{(.)})$.
The fragmentation is homogeneous in time if $P_{t,t'}$ depends
only on $t'-t$.
\end{itemize}
A Markov process $(\Gamma(t),t\ge 0)$ on $\sr$ is called an
exchangeable fragmentation process of  compositions if, for all
$n\in \mcn$, the process $(\Gamma_{[n]}(t),t\ge 0)$   is an exchangeable
fragmentation process of  compositions on $\sr_n$.
\end{defi}

In the sequel, a $c$-fragmentation will denote an exchangeable fragmentation
process on compositions.

\subsection{Interval fragmentation}
In this section we recall the definition of a
homogeneous\footnote{In \cite{Bertoin02}, Bertoin defines more
generally self-similar interval fragmentations with index $\alpha$.
Here, the term homogeneous means that we only
consider the case $\alpha=0$. } interval fragmentation ~\cite{Bertoin02}.\\
We consider a family of probability measures $(q_{t,s},t\ge 0,s>t)$
on $\ur$.
 For all interval $I=]a,b[\subset]0,1[$, we define the affine transformation
$g_I:]0,1[\rightarrow I$ given by $g_I(x)=a+x(b-a)$. We still denote
$g_I$ the induced map on $\ur$, so, for $V\in \ur$, $g_I(V)$ is an
open subset of $I$. We define then  $q_{t,s}^I$ as the image of
$q_{t,s}$ by $g_I$. Hence $q_{t,s}$ is a probability measure on the
open subset of $I$. Finally, for $W\in \ur$ with interval
decomposition $(I_i,i\in\mcn)$, $q_{t,s}^W$ is the distribution of
$\cup X_i$ where the $X_i$ are independent random variables with
respective law $q_{t,s}^{I_i}$.

\begin{defi}
A process $(U(t),t\ge 0)$ on $\ur$ is called a  homogeneous
interval fragmentation if it is a Markov process which fulfills
the following properties :
\begin{itemize}
\item $U$ is continuous in probability and $U(0)=\un$ a.s. \item
$U$ is nested i.e. for all $s>t$ we have $U(s)\subset U(t)$. \item
There exists a family $(q_{t,s},t\ge 0,s>t)$ of probability measure
on $\ur$ such that :
$$\forall t\ge 0,\; \forall s>t,\; \forall A\subset \ur,\quad  \mcp(U(s)\in
A|\; U(t))=q_{t,s}^{U(t)}(A).$$
\end{itemize}
In the following, we abbreviate an interval fragmentation process as  an
$i$-fragmentation.
\end{defi}

We remark that if we take the decreasing sequence of the sizes of
the
 interval components of an $i$-fragmentation, we obtain a mass-fragmentation, denoted here a $m$-fragmentation
  (see \cite{CoursBertoin03} for definition of $m$-fragmentation).
But, with the $m$-fragmentation, we loose the genealogical aspect
present in the $i$-fragmentation.

\subsection{Link between $i$-fragmentation and $c$-fragmentation }

From this point of the paper and until Section \ref{inhomogene}, the
fragmentation processes will always be homogeneous in time, i.e.
$q_{t,s}$ depends only on $s-t$, hence we will just write  $q_{t-s}$
to denote $q_{t,s}$.

\begin{theo}\label{ifragcgrag}
There is a one to one correspondence between laws of $i$-fragmentations
and laws of $c$-fragmentations. More precisely : \begin{itemize}
\item let $(U(t),t\ge 0)$ be an $i$-fragmentation. Let $(V_i)_{i\ge
0}$ be a sequence of independent random variables uniformly
distributed on ]0,1[. Using the same process as in Definition
\ref{paintbox} with $U(t)$ and $(V_i)_{i\ge 1}$, we define a process
$(\Gamma(t),t\ge 0)$ on $\sr$. Then $(\Gamma(t),t\ge 0)$ is a
$c$-fragmentation and we have $U_{\Gamma(t)}=U(t)$  a.s. for each
$t\ge 0$. \item Let $(\Gamma(t),t\ge 0)$ be a $c$-fragmentation.
Then $(U_{\Gamma(t)},t\ge 0)$ is an $i$-fragmentation.
\end{itemize}

\end{theo}

Proof. We begin by proving the first point. We have by Theorem
\ref{gnedinouv}, $U_{\Gamma(t)}=U(t)$ a.s. for each $t\ge 0$. Let us
fix $n\in \mcn$ and    $t\ge 0$.  We are going to prove that, for
$s>t$, the conditional law of $\Gamma_{[n]}(s)$ given
$\Gamma_{[n]}(t)=(\gamma_1,\ldots,\gamma_k)$ is the law of
$FRAG(\Gamma_{[n]}(t),\gamma^{(\cdot)})$, where $\gamma^{(\cdot)}$
is
 a  sequence of iid exchangeable compositions with law $\Gamma_{[n]}(s-t)$. Since
 $(U(s),s\ge 0)$ is a fragmentation process, we have $U(t+s)\subset U(t)$. By construction of $\Gamma_{[n]}(t)$, it is then clear
 that $\Gamma_{[n]}(t+s)$ is a finer composition than $\Gamma_{[n]}(t)$. Hence each singleton of $\Gamma_{[n]}(t)$ remains a singleton of $\Gamma_{[n]}(t+s)$.
  So we can assume that $\Gamma_{[n]}(t)$ has no singleton. For $1\le i\le k$, fix $l\in \gamma_i$ and  define
  $$a_i=\sup\{a\le V_l,a\notin U(t)\}, \qquad b_i=\inf\{b\ge V_l,b\notin U(t)\}.$$
  Notice that $a_i$ and $b_i$ do not depend on the choice of $l\in \gamma_i$. Furthermore, since $\Gamma_{[n]}(t)$ has no singleton,
  we have $a_i<b_i$ almost surely.
   We define also
  $$Y^i_j=\left(\frac{V_j-a_i}{b_i-a_i}\right)_{j\in \gamma_i,1\le i \le k}.$$
By construction of $\Gamma_{[n]}(t)$, the random variables $(Y^i_j)_{j\in \gamma_i,1\le i \le k}$ are independent and uniformly distributed on $]0,1[$.
Besides, $(]a_i,b_i[)_{1\le i\le k}$
 are $k$ distinct interval components of $U(t)$. Since $U(t)$ is a fragmentation process, the processes
 $$\left(U^i(s)=\frac{1}{b_i-a_i}(U_{]a_i,b_i[}(s)-a_i),s\ge t\right)_{1\le i\le k}$$ are $k$ independent $i$-fragmentations with law $(U(s-t),s\ge t)$.
 Let $\gamma^{(i)}(s)$ be the composition of $\gamma_i$ obtained from $U^i(s)$ and $(Y^i_j)_{j\in \gamma_i}$ using Definition \ref{paintbox}.
Hence, $\gamma^{(i)}(s)$ has the law of $\Gamma_{\gamma_i}(s-t)$ and the processes $(\gamma^{(i)}(s),s\ge t)_{1\le i\le k}$ are independent.
Furthermore, by construction we have $\Gamma_{[n]}(t+s)=FRAG(\Gamma_{[n]}(t),\gamma^{(\cdot)}(s))$. Hence, $(\Gamma_{[n]}(t),t\ge 0)$ has the expected
 semi group of transition.

\vspace*{0.3cm}

  Let us now prove the second point. In the following, we will
write $U_t$ to denote $U_{\Gamma(t)}$. First, we prove that for all
$s>t$, $U_s\subset U_t$. Fix $x \notin U_t$, we shall prove  $x
\notin U_s$. We have $\chi_{ U_t}(x)=\min\{|x-y|,y\in U_t^{c}\}=0$.
Let $U^n_t$ be the open subset of $]0,1[$ corresponding to
$\Gamma_{[n]}(t)$ as in Theorem \ref{gnedinouv}. So we have
$\lim_{n\rightarrow \infty} d(U^n_t,U_t)=0$. Fix $\varepsilon >0$.
Hence, there exists $N\in\mcn$ such that, for all $n\ge N$, $\chi_{
U^n_t}(x)\le \varepsilon$. This implies that :

$$\forall n\ge N, \exists y_n\notin U_t^n \mbox{ such that } |y_n-x|\le \varepsilon.$$
Besides, as  $(\Gamma(t),t\ge 0)$ is a fragmentation, we have for
all $n\in \mcn$, $U_s^n\subset U_t^n$. Hence, we have also
$$\forall n\ge N,  y_n\notin U_s^n,$$
and so $\chi_{ U^n_s}(x)\le \varepsilon$ for all $n\ge N$. We deduce
that $\chi_{ U_s}(x)=0$ i.e.
 $x \notin U_s$.

We now have to prove the branching property. Fix $t>0$. We
consider the decomposition of $U_t$  in disjoint intervals :
$$U_t=\coprod_{k\in \mcn} I_k(t).$$
 Set $F_k(s)=U_{t+s}\cap I_k(t)$.
We want to prove that, given $U_t$  :
\begin{itemize}
\item  $\forall l\in \mcn, \forall m_1,\ldots, m_l$ distinct,
$F_{m_1},\ldots,F_{m_l}$ are independent processes. \item $F_k$
has the following law :
$$\forall A \mbox{ open subset of }]a,b[,\; \mcp((F_k(s),s\ge 0) \in A \;| I_k(t)=]a,b[)=\mcp((U_{s},s\ge 0)\in (b-a)A+a).$$
\end{itemize}

For all $k\in \mcn$, there exists $i_k\in \mcn$ such that, if
$J_{i_k}^n(t)$ denotes the interval component of $U_t^n$ containing
the integer $i_k$, then $J_{i_k}^n(t) \overset{n\rightarrow
\infty}{\longrightarrow} I_k(t)$. Let $B_k$ be the block of
$\Gamma(t)$ containing $i_k$. As $B_k$ has a positive asymptotic
frequency,
 it is isomorphic to $\mcn$. Let $f$ be the increasing bijection from the set of element of $B_k$ to $\mcn$.
  Let us re-label the elements of $B_k$ by their image by $f$. The
  process
$(U_{\Gamma_{B_k}(t+s)},s\ge 0)$ has then the same law as
$(U_{s},s\ge 0)$ and  is independent of the rest of the
fragmentation. Besides, given $I_k(t)=]a,b[$,
$F_k(s)=a+(b-a)U_{\Gamma_{B_k}(t+s)}$, so the two points above are
proved. $\Box$

Hence, this result completes an analogous result due to Berestycki
\cite{Berestycki02} in the case of $m$-fragmentations and
$p$-fragmentation (i.e. fragmentations of exchangeable partitions).
We can again draw a diagram to represent the link between the four
kinds of fragmentation :

$$
\begin{array}{ccc}
 \big(\sr,(\Gamma(t),t\ge 0)\big)&\xleftrightarrow{\mbox{\tiny{ Theorem \ref{ifragcgrag} }}}& \big(\ur,(U_{\Gamma(t)},t\ge 0)\big)\\
\mbox{\tiny{$\wp_1$}}\Big{\downarrow } & & \mbox{\tiny{$\wp_2$}}\Big{\downarrow }  \\
 \big(\pinf,(\Pi(t),t\ge 0)\big)
&\xleftrightarrow{\hspace*{0.2cm}\mbox{\tiny{Berestycki}\quad}} &
\big(\srb,(\freq{U_{\Gamma(t)}},t\ge 0)\big).
 \end{array}
$$

\section{Some general properties}
In this section, we gather general properties of $i$ and
$c$-fragmentations. Since the proof of these results are simple
variations of those in the case of $m$ and $p$-fragmentations
\cite{CoursBertoin03}, we will be a bit sketchy.

\subsection{Measure of a fragmentation process}

Let $(\Gamma(t),t\ge 0)$ be a $c$-fragmentation.
 As in the case of $p$-fragmentation  \cite{CoursBertoin03}, for $n\in \mcn $ and $\gamma \in \sr_n^{*}$, we define a jump rate
from $\un_{n}$ to $\gamma$ :
$$q_{\gamma}=\lim_{s\rightarrow 0}\frac{1}{s}\mcp\left(\Gamma_{[n]}\left(s\right)=\gamma
\right).$$ With the same arguments as in the case of
$p$-fragmentation, we can also prove that the family
$(q_{\gamma},\gamma \in \sr_n^{*}, n\in \mcn)$ characterizes the law
of the fragmentation (you just have to use that distinct blocks
evolve independently and with the same law). Furthermore,
 observing that we have
$$\forall n<m, \quad \forall \gamma'\in \sr_n^{*},\quad q_{\gamma'}=\sum_{\gamma\in \sr_m,
\gamma_{[n]}=\gamma'}q_{\gamma},$$ and that
$$\forall n\in \mcn, \quad \forall \sigma \in \mbox{\LARGE{$\sigma$}\normalsize$_n$}, \quad \forall \gamma \in \sr_n^{*},\quad
q_{\gamma}=q_{\sigma(\gamma)},$$
we deduce  that there exists a unique exchangeable measure $\mu$ on
$\sr$ such that $\mu(\mathbf{1})=0$ and
$\mu(\mathcal{Q}_{\infty,\gamma})=q_\gamma$ for all $\gamma \in
\sr_n^*$ and $n \in \mcn$, where
$\mathcal{Q}_{\infty,\gamma}=\{\gamma' \in \mathcal{\sr},
\gamma'_{[n]}=\gamma\}$. Furthermore, the measure $\mu$ characterizes
the law of the fragmentation.

We remark also that if $\mu$ is the measure of a fragmentation
process, we have for all $n \in \mcn$,
$$\mu(\{\gamma \in \sr, \gamma_{[n]}\neq
\un_n \})=\sum_{\gamma \in \sr_n^*}{q_{\gamma}}<\infty.$$ So we can
apply Theorem \ref{disc} to $\mu$ and we deduce  the following
corollary :

\vspace*{0.3cm}

\begin{cor}
Let $\mu$ be the measure of a $c$-fragmentation. Then there exist a
dislocation measure $\nu$ and two nonnegative numbers $c_l$ and
$c_r$ such that :
\begin{itemize}
\item $\mu \un_{\{ U_{\gamma}\neq\un\}}=P^{\nu}$. \item $\mu
\un_{\{ U_{\gamma}=\un\}}=c_l\epsilon+c_r\epsilon'$.
\end{itemize}
\end{cor}

With a slight abuse of notation, we will write sometimes in the sequel that
 $\mu=(\nu,c_l,c_r)$ when $\mu=P^{\nu}+c_l\epsilon+c_r\epsilon'$.

\subsection{The Poissonian construction}

Let us recall that we define in Section  \ref{mesure} a
fragmentation measure as a measure $\mu$ on $\sr$ such that :
\begin{itemize}
\item $\mu$ is  exchangeable.
 \item $\mu(\un_{\mcn})=0.$
 \item For all $n \in \mcn$,  $\mu(\{\gamma \in \sr,
\gamma_{[n]}\neq \un_n \})<\infty$.
\end{itemize}

Notice that if $\mu$ is the measure of a $c$-fragmentation, then
$\mu$ is a fragmentation measure. Conversely, we now prove that, if  we consider a fragmentation measure $\mu$,  we can
construct a $c$-fragmentation with measure
  $\mu$.

We consider a Poisson measure $M$ on
$\mathbb{R}_{+}\times\sr\times\mcn$ with intensity $dt\otimes\mu\otimes \sharp$,  where
$\sharp$ is the counting measure on $\mcn$ .
Let  $M^n$ be the restriction of $M$ to
$\mathbb{R}_{+}\times\sr_n^*\times\partn$. The  intensity measure is
then finite on the interval $[0,t]$,
so we can  order the atoms of $M^n$ according to their first coordinate. \\
For $n\in \mcn$, $(\gamma,k)\in \sr\times \mcn$, let
$\Delta^{(.)}_n(\gamma,k)$ be the composition sequence of $\sr_n$
defined by :
$$\Delta^{(i)}_n(\gamma,k)=\un_{n} \quad \mbox{ if } i\neq k   \hspace*{0.8cm}\mbox{ and }
   \hspace*{0.8cm} \Delta^{(k)}_n(\gamma,k)=\gamma_{[n]}.$$
We construct then a process $(\Gamma_{[n]}(t),t\ge 0)$
on $\sr_n$ in the following way :\\
$\Gamma_{[n]}(0)=\un_n$.\\
$(\Gamma_{[n]}(t),t\ge 0)$ is a pure jump process which jumps at times
 when an atom of $M^n$ appears. More precisely, if $(s,\gamma,k)$ is an atom  of $M^n$,
 set
$\Gamma_{[n]}(s)=FRAG(\Gamma_{[n]}(s^-),\Delta^{(.)}_n(\gamma,k))$. \\
We can check that this construction is compatible with the
restriction ; hence, this defines a process $(\Gamma(t),t\ge 0)$ on $\sr$.

\begin{prop}
Let $\mu$ be a fragmentation measure. The construction above  of a
process on  compositions from a Poisson point process on $\mcr_+\times
\sr\times\mcn$ with intensity $dt\otimes\mu\otimes \sharp$, where
$\sharp$ is the counting measure on $\mcn$ , yields a
$c$-fragmentation with measure $\mu$.
\end{prop}
The proof is an easy adaptation of the Poissonian construction of $p$-fragmentations (cf. \cite{CoursBertoin03}).
As the sequence $\Delta^{(.)}_n(\gamma,k)$ is
doubly exchangeable, we also have that $\Gamma_{[n]}(t)$ is an
exchangeable composition for each $t\ge 0$. Looking as the rate jump
of the process $\Gamma_{[n]}(t)$, it is then easy to check that the
constructed process is a $c$-fragmentation with measure $\mu$.$\Box$

\vspace*{0.3cm}

A Poissonian construction of an $i$-fragmentation with no erosion is
also possible with a Poisson measure on
$\mathbb{R}_{+}\times\ur\times\mcn$ with intensity
$dt\otimes\nu\otimes \sharp$. The proof of this result  is not as
simple as for compositions because we can not restrict to a discrete
case as done above. In fact, to prove this proposition, we must take
the image of the Poisson measure M above by an appropriate
application. For more details refer to Berestycki
\cite{Berestycki02} who have already proved this result for
$m$-fragmentation and the same approach works in our case.

\vspace*{0.5cm}

To conclude this section, let us notice how the two erosion coefficients
affect the fragmentation.
Let $(U(t),t\ge 0)$ be an $i$-fragmentation
with parameter $(0,c_l,c_r)$. Set $c=c_l+c_r$. We have :
$$U(t)= \left]\frac{c_l}{c}(1-e^{-tc}),1-\frac{c_r}{c}(1-e^{-tc})\right[ \mbox{ a.s.}$$
Indeed, consider a  $c$-fragmentation $(\Gamma(t),t\ge 0)$ such
that $u_{\Gamma(t)}=U(t)$ a.s. We define
$\mu_{c_l,c_r}=c_l\epsilon + c_r\epsilon'$. Hence
$(\Gamma(t),t\ge 0)$ is a fragmentation with measure $\mu_{c_l,c_r}$.
Recall that the process $(\Gamma(t),t\ge 0)$ can be constructed from
a Poisson measure on $\mathbb{R}_{+}\times \sr \times \mcn$ with
intensity $dt\otimes \mu_{c_l,c_r} \otimes \sharp$. By the form of
$\mu_{c_l,c_r}$, we remark then that, for all $t\ge0$, $\Gamma(t)$
have only one block non reduced to a singleton. Furthermore,
 for all $n\in \mcn$, the integer $n$
is a singleton at time $t$ with probability $1-e^{-tc}$, and,
given $n$ is a singleton of $\Gamma(t)$, $\{n\}$ is before the
infinite block of $\Gamma(t)$ with probability $c_l/c$ and after with
probability $c_r/c$. By the law of large number, we deduce that
the proportion of singletons before the infinite block of $\Gamma(t)$
is almost surely $\frac{c_l}{c}(1- e^{-tc})$ and the
 proportion of singletons after the infinite block of $\Gamma(t)$ is almost surely $\frac{c_r}{c}(1-e^{-tc})$.

\begin{rem}
Berestycki ~\cite{Berestycki02} has proved a similar result for
the $m$-fragmentation. He also proved  that if $(F(t),t\ge 0)$ is
a $m$-fragmentation with parameter $(\nu,0)$, then
$\tilde{F}(t)=e^{-ct}F(t)$ is a $m$-fragmentation with parameter
$(\nu,c)$. But, we can not generalize this result for $i$-fragmentation
 because the proportion of singleton between two
successive component intervals of the fragmentation depends on the
time where the two component intervals split.
\end{rem}

\subsection{Projection from $\ur$ to $\srb$}

We know that if $(U(t),t\ge 0)$ is an $i$-fragmentation, then its
projection on $\srb$, $(\freq{U(t)},t\ge 0)$ is a $m$-fragmentation.
More precisely, we can express the characteristics of the
$m$-fragmentation from the characteristics of the $i$-fragmentation.

\begin{prop}\label{massinter}
The ranked sequence of the length of an
$i$-fragmentation with measure $(\nu,c_l,c_r)$ is a $m$-fragmentation with parameter
$(\tilde{\nu},c_l+c_r)$ where $\tilde{\nu}$ is the
image of $\nu$ by the application $U\rightarrow
|U|^{\downarrow}$.
\end{prop}

Proof. Let $(\Gamma(t),t\ge 0)$ be a $c$-fragmentation with measure $\mu=(\nu,c_l,c_r)$. Let $(\Pi(t),t\ge 0)$ be its image by $\wp_1$.
The process  $(\Pi(t),t\ge 0)$
is then a $p$-fragmentation. Set $n\in\mcn$ and $\pi\in \pn^*$.
We have
\begin{eqnarray*}
q_{\pi}&=&\lim_{s\rightarrow 0}\frac{1}{s}\mcp(\Pi_{[n]}(s))=\pi)\\
&=&\lim_{s\rightarrow 0}\frac{1}{s}\mcp\left(\Gamma_{[n]}(s))\in \wp^{-1}(\pi)\right)\\
&=&\tilde{\mu}(\pi),
\end{eqnarray*}
where $\tilde{\mu}$ is the image of $\mu$ by $\wp_1$. Besides we have already prove that $\tilde{\mu}=(\tilde{\nu},c_l+c_r)$.
We consider now the $i$-fragmentation $(U_{\Gamma(t)},t\ge 0)$ with measure $(\nu,c_l,c_r)$. We get
that the process $(\freq{U_{\Gamma(t)}},t\ge 0)$ is a.s. equal to the $m$-fragmentation $(\freq{\Pi(t)},t\ge 0)$ which fragmentation
measure is $(\tilde{\nu},c_l+c_r)$. $\Box$

\vspace*{0.7cm}

According to Proposition \ref{massinter} and using the theory of
$m$-fragmentation (see \cite{CoursBertoin03}), we deduce then the
following results :
\begin{itemize}
\item Let $(\Gamma(t),t\ge 0)$ be a $c$-fragmentation  with
parameter $(\nu,c_l,c_r)$. We denote by $B_1$ the block of $\Gamma(t)$ containing the integer $1$.
 Set $\sigma(t)=-\ln |B_1(t)|$. Then $(\sigma(t),t\ge 0)$ is a subordinator. If we denote $\zeta=\sup\{t>0,
\sigma_t<\infty\}$, then there exists a non-negative
function $\phi$ such that
$$\forall q,t\ge 0,\;\, \mce[\exp(-q\sigma_t),\zeta>t]=\exp(-t\phi(q)du).$$
We call $\phi$  the Laplace exponent of  $\sigma$ and we have :
$$\phi(q)=(c_l+c_r)(q+1)+\int_{\ur}{(1-\sum_{i=1}^{\infty}{|U_i|^{q+1}})\nu(dU)},$$
where  $(|U_i|)_{i\ge 0}$ is the sequence of the lengths of the
component intervals of
 $U$.


\item An $(\nu,c_r,c_l)$ $i$-fragmentation
$(U(t),t\ge 0)$ is proper (i.e. for each $t$, $U(t)$ has almost
surely a Lebesgue measure equal to $1$) iff
$$ c_l=c_r=0 \mbox{ and } \nu\left(\sum_{i}{s_i}<1\right)=0.$$
\end{itemize}

\subsection{Extension to the time-inhomogeneous case}\label{inhomogene}
We now briefly expose how the results of the preceding sections can be transposed in the case of time-inhomogeneous
fragmentation. We will not always detail the proof since their are very similar as in the homogeneous case.
In the sequel, we shall focus on c-fragmentation $(\Gamma(t),t\ge 0)$
fulfilling the following properties :
\begin{itemize}
\item for all $n \in \mcn$, let $\tau_n$ be the time of the first
jump of  $\Gamma_{[n]}$ and $\lambda_n$ be its law. Then $\lambda_n$ is
absolutely  continuous with respect to Lebesgue measure with
continuous and strictly positive density. \item for all $\gamma \in
\sr_n^*,\hspace*{0.3cm}
h_{\gamma}^n(t)=\mcp(\Gamma_{[n]}(t)=\gamma\;|\;\tau_n=t)$ is a continuous
function of t.
\end{itemize}

Remark that a time homogeneous fragmentation always fulfills this
two points. Indeed, in that case, $\lambda_n$ is an exponential
random variable and the function $h_{\gamma}^n(t)$ does not depend on
$t$.
 As in the case of fragmentation of exchangeable partitions \cite{Basdevant05},
 for $n\in \mcn $ and $\gamma \in \sr_n^{*}$, we can define an instantaneous rate of jump
from $\un_n$ to $\gamma$ :
$$q_{\gamma,t}=\lim_{s\rightarrow 0}\frac{1}{s}\mcp\left(\Gamma_{[n]}\left(\tau_n\right)=\gamma
\; \& \;\tau_n \in [t,t+s]\;|\; \tau_n\ge t\right).$$ With the same
arguments as in the case of  fragmentation of exchangeable
partitions \cite{Basdevant05}, we can also prove that, for each
$t>0$,  there exists a unique exchangeable measure $\mu_t$ on $\sr$
such that  $\mu_t(\mathbf{1})=0$ and $\mu_t(\mathcal{Q}_{\infty,\gamma})=q_{\gamma,t}$
for all $\gamma \in \sr_n\backslash \{\mathbf{1}\}$ and $n \in \mcn$,
where $\mathcal{Q}_{\infty,\gamma}=\{\gamma' \in \mathcal{\sr}, \gamma'_{n}=\gamma\}$.
Furthermore, the family of measure $(\mu_t,t\ge 0)$ characterizes
the law of the fragmentation.

We remark also that if $(\mu_t,t\ge 0)$ is the family of measure of
a fragmentation process, we have for all $n \in \mcn$,
$$\mu_t(\{\gamma \in \sr, \gamma_{[n]}\neq
\mathbf{1}_n \})=\sum_{\gamma \in \sr_n^*}{q_{\gamma,t}}<\infty \mbox{ and }
\int_{0}^{t}{\mu_u(\{\gamma \in \sr, \gamma_{[n]}\neq \mathbf{1}_n
\}}=-\ln(\lambda_n(]t,\infty[))<\infty.$$ So we can apply Theorem
\ref{disc} to $\mu_t$ and we deduce  the following proposition :

\begin{cor}
Let $(\mu_t,t\ge 0)$ be the family of measure of a $c$-fragmentation.
Then there exists a family of dislocation measures
$(\nu_t,t\ge 0)$ and two families of nonnegative numbers $(c_{l,t},t\ge
0)$, $(c_{r,t},t\ge 0)$ such that :
\begin{itemize}
\item $\mu_t \un_{\{ U_{\pi}\neq\un\}}=P^{\nu_t}$. \item $\mu_t
\un_{\{ U_{\pi}=\un\}}=c_{r,t}\epsilon+c_{r,t}\epsilon'$.
\end{itemize}
Besides we have for all $T\ge 0$,
\begin{equation*}
\int_{0}^{T}{\int_{\ur}{\left(1-s_1\right)\nu_{t}\left(dU\right)dt}<\infty}
\mbox{ and }\int_{0}^{T}{(c_{l,t}+c_{r,t} )dt}<\infty.
\end{equation*}
\end{cor}

 The first part of the proposition comes from Theorem \ref{disc}.  For the second part, use that

$$\int_{\ur}{\left(1-s_1\right)\nu_t\left(dU\right)}\le \mu_t\left(\{\pi \in \pinf, \pi_{|2}\neq \mathbf{1}
\}\right).$$
 For the upper bound concerning  the erosion coefficients, we remark that  :
$$c_t+c'_t=\mu_t\left(\{1\},\mcn\setminus\{1\}\right)+
\mu_t\left(\mcn\setminus\{1\},\{1\}\right).\Box$$

In the same way as for homogeneous fragmentation, we define  a fragmentation measure family as  a family $(\mu_t,t\ge 0)$  of
exchangeable measures  on $\sr$ such that, for
each $t\in [0,\infty[$, we have :
\begin{itemize}
\item $\mu_t(\un_{\mcn})=0$. \item $\forall n \in \mcn \; \mu_t(\{\gamma \in
\sr, \gamma_{[n]}\neq \un_n \})<\infty$ and
$\int_{0}^{t}{\mu_u(\{\gamma \in \sr, \gamma_{[n]}\neq \un_n
\}})du<\infty$. \item $\forall n \in \mcn, \; \forall A\subset
\sr_n^*,\; \mu_t(A)$ is a continuous function of $t$.
\end{itemize}

\begin{prop}
Let $(\mu_t,t\ge 0)$ be a fragmentation measure family. A $c$-fragmentation with fragmentation measure
$(\mu_t,t\ge 0)$  can be constructed from a Poisson point process on  $\mcr_+\times
\sr\times\mcn$ with intensity $dt\otimes\mu_t\otimes \sharp$, where
$\sharp$ is the counting measure on $\mcn$ in the same way as for time-homogeneous fragmentation.
\end{prop}

It is very easy to check that the proof of the homogeneous case applies here too.

Of course,
a Poissonian construction of a time-inhomogeneous $i$-fragmentation with no
erosion is also possible with a Poisson measure on
$\mathbb{R}_{+}\times\ur\times\mcn$ with intensity
$dt\otimes\nu_t\otimes \sharp$.

Concerning the law of the tagged fragment, if you define $\sigma(t)=-\ln |B_1(t)|$, with $B_1$ the block containing
the integer $1$, we have now that $\sigma(t)$  is a process with
independent increments. And so, if we denote $\zeta=\sup\{t>0,
\sigma_t<\infty\}$, then there exists a family of non-negative
functions $(\phi_t,t\ge 0)$ such that
$$\forall q,t\ge 0,\;\, \mce[\exp(-q\sigma_t),\zeta>t]=\exp(-\int_{0}^{t}\phi_u(q)du).$$
We call $\phi_t$  the instantaneous Laplace exponent of  $\sigma$
at time $t$ and we have :
$$\phi_t(q)=(c_{l,t}+c_{r,t})(q+1)+\int_{\ur}{(1-\sum_{i=1}^{\infty}{|U_i|^{q+1}})\nu_t(dU)},$$
where  $(|U_i|)_{i\ge 0}$ is the sequence of the lengths of the
component intervals of
 $U$.
 Furthermore,
 an $(\nu_t,c_t,c_t')_{t\ge 0}$ $i$-fragmentation
$(U(t),t\ge 0)$ is proper iff :
$$\forall t>0, \quad c_{l,t}=c_{r,t}=0 \quad \mbox{ and }\quad \nu_t(\sum_{i}{s_i}<1)=0).$$

Finally, we can also compute the law of an $(0,c_{l,t},c_{r,t})_{t\ge0}$ $i$-fragmentation.
 After some calculus, we obtain that
we have :
$$U(t)=\Big]\int_{0}^{t}{c_{l,u}\exp(-C_u)du},\; 1-\int_{0}^{t}{c_{r,u}\exp(-C_u)du}\Big[ \mbox{ a.s.}$$
with $C_u=\int_{0}^{u}(c_{l,v}+c_{r,v}) dv$.

\subsection{Extension to the self-similar case}
A notion of self similar fragmentations has been also introduced  \cite{Bertoin02}.
We recall here the definition of a self similar $p$-fragmentation, the reader can easily adapt this definition to the three
other cases of fragmentation.

\begin{defi}
Let $\Pi=(\Pi(t),t\ge 0)$ be an exchangeable process on $\pinf$. We call $\Pi$ a self similar $p$-fragmentation with index $\alpha\in \mcr$ if
\begin{itemize}
\item $\Pi(0)=1_{\mcn}$ a.s.
\item $\Pi$ is continuous in probability
\item For every $t\ge 0$, let $\Pi(t)=(\Pi_1,\Pi_2,\ldots)$  and denote by
$|\Pi_i|$ the asymptotic frequency of the block $\Pi_i$. Then for every $s>0$,
the conditional distribution of $\Pi(t+s)$ given $\Pi(t)$ is the law of the random
partition whose blocks are those of the partitions $\Pi^{(i)}(s_i)\cap \Pi_i$ for $i\in \mcn$, where $\Pi^{(1)},\ldots$ is a sequence of
independent copies of $\Pi$ and $s_i=s|\Pi_i|^\alpha$.
\end{itemize}
\end{defi}

Notice that an homogeneous $p$-fragmentation corresponds to the case $\alpha=0$.

We have still the same correspondence between the four types of fragmentation. In fact, a self similar fragmentation can be constructed from
a homogeneous fragmentation with a time change :

 \begin{prop}\label{change}\cite{Bertoin02} Let $(U(t),t\ge 0)$ be an
 homogeneous interval fragmentation with measure $\nu$. For $x\in
 ]0,1[$, we denote by $I_x(t)$ the interval component of $U(t)$
 containing $x$. We define  $$T_t^{\alpha}(x)=\inf\{u\ge 0,
 \int_{0}^{u}|I_x(r)|^{-\alpha}dr>t\}\mbox{ and } U^{\alpha}(t)=U(T_t^{\alpha})=\bigcup
 I_x(T_t^{\alpha}(x)).$$
Then $(U^\alpha(t),t\ge 0)$ is a self similar interval
fragmentation with index $\alpha$.
 \end{prop}

A self similar $i$-fragmentation (or $c$-fragmentation) is then characterized by a quadruplet $(\nu,c_l,c_r,\alpha)$
where $\nu$ is a dislocation measure on $\ur$, $c_l$ and $c_r$ are two nonnegative numbers and $\alpha\in \mcr$ is the index of self similarity.

\section{Hausdorff dimension of an interval fragmentation}\label{Haussec}

Let $(U(t),t\ge 0)$ be a self similar $i$-fragmentation with index
$\alpha>0$. Let $K(t)=[0,1]\backslash U(t)$. The set $K(t)$ is a
closed set, and if the fragmentation is proper (i.e. the
fragmentation has with no erosion and its fragmentation measure
verifies $\nu(\sum_i \freq{U_i}<1)=0$), its Lebesgue measure is
equal to 0. Hence, to evaluate the size of $F(t)$, we shall compute
its Hausdorff measure. Here, we will just examine time-homogeneous
fragmentation. First we recall the definition of the Hausdorff
dimension of a subset of ]0,1[.

\begin{defi}\cite{Falconer86}\label{hauss}
Let $A\in]0,1[$.  Let $d\ge 0$ and $r>0$. We set
$$J_d^r(A)=\inf\left\{\sum_{i=1}^{\infty}|b_i-a_i|^d, A\subset \bigcup_{i=1}^{\infty}[a_i,b_i],|b_i-a_i|\le r
\right\} \mbox{ and }H_d(A)=\lim_{r\rightarrow 0^+}J_d^r(A),$$
(this limit exists since $J_d^r(A)$ decreases with $r$). $H_d(A)$
is the $d$-Hausdorff measure of $A$. Furthermore, there exists a
unique number $D$ such that
$$\forall d>D, H_d(A)=0 \mbox{ and } \forall d<D, H_d(A)=\infty.$$
This number is the Hausdorff dimension of $A$ and is denoted by
$\dim_{\mathcal{H}}(A)$.
\end{defi}


We will now calculate the Hausdorff dimension of the complement of
a time-homogeneous $i$-fragmentation in the case where the measure of
fragmentation fulfils some conditions.

\begin{hypo}\label{hyp}
Let  $\nu$ be a dislocation measure. We  assume that $\nu$ fulfills
the following conditions  :
\begin{itemize}
\item[(H1)] $\nu$ is conservative i.e. $\nu(\sum_i \freq{U_i}<1)=0$.
 \item[(H2)] There exists an integer $k$ such that $\nu(\freq{U_k}>0)=0$, i.e. $\nu$ is carried by the open sets with at most $k-1$ interval components.
 \item[(H3)] Let $h(\varepsilon)=\int_{\ur}(Card\{i, |U_i|\ge
 \varepsilon\}-1) \nu(dU)$. Then $h$ is regularly
 varying  with index $-\beta$ as $\varepsilon \rightarrow 0+$.
\item[(H4)] Let $g$ be the left extremity of the largest interval component of a generic open set and $d$
 the right extremity. Then  as $ \varepsilon\rightarrow 0+$,
 we have either $\liminf \frac{\nu(g\ge \varepsilon)}{ \nu(d\le 1-\varepsilon)}>0$ or
$\limsup \frac{\nu(g\ge \varepsilon)}{ \nu(d\le
1-\varepsilon)}<\infty$.
\end{itemize}

\end{hypo}

We can now state the theorem :
\begin{theo}\label{theohaus} Let $\nu$ be a dislocation measure fulfilling  Hypothesis \ref{hyp}.
Let $(U(t),t\ge 0)$ be an $i$-fragmentation with characteristics
$(\nu,0,0)$   and index of self-similarity $\alpha$ strictly
positive. Let $K(t)=[0,1]\backslash U(t)$. Then the Hausdorff
dimension of $K(t)$ is $\beta$ for all $t>0$ simultaneously, a.s.
\end{theo}

In fact, if the index of self-similarity is zero, the lower bound of
the Hausdorff dimension still holds. Besides, Hypothesis $(H4)$ is
only needed to prove the lower bound and allows a large class of
dislocation measure such as symmetric measures or, at the opposite,
measures for which the largest fragment is always on the same side.

\vspace*{0.5cm}

Proof. We will first prove the upper bound. Let us recall a lemma
 proved by Bertoin in \cite{Bertoin04} for $m$-fragmentation
processes whose dislocation measure fulfills Hypothesis \ref{hyp}.

\begin{duge}\cite{Bertoin04}\label{dugehaus}
Let $(U(t),t\ge 0)$ be  a self-similar $(\nu,0,0,\alpha)$
$i$-fragmentation with index of self similarity strictly positive
and whose dislocation measure fulfills  (H1), (H2), (H3). Let
$(X(t)=(X_i(t))_{i\ge 1},t\ge 0)$ the associated $m$-fragmentation.
Let $N(\varepsilon,t)=Card\{i\ge 1, X_i(t)\ge \varepsilon\}$ and
$M(\varepsilon,t)=\sum_i{X_i(t)\un_{\{X_i\le \varepsilon\}}}$. Then
$\lim_{\varepsilon\rightarrow
0+}\frac{N(\varepsilon,t)}{h(\varepsilon)} \mbox{ and }
\lim_{\varepsilon\rightarrow 0+}\frac{M(\varepsilon,t)}{\varepsilon
h(\varepsilon)}$ exist and are strictly positive and finite.
\end{duge}

Let us now fix $d\in ]0,1[$ and look for a upper bound  of the
$d$-Hausdorff measure of $K(t)$. Let
$I_\varepsilon=]0,1[\backslash\{\mbox{ interval components of
}U(t)\mbox{ which size is larger than }\varepsilon\}$. So we have
$K(t)\subset I_{\varepsilon}$ and
$|I_{\varepsilon}|=M(\varepsilon,t)$ since $\nu$ is conservative.
Furthermore, $I_{\varepsilon}$ has at most $N(\varepsilon,t)+1$
interval components. Using notation of Definition \ref{hauss}, we
get :
$$J_d^{\varepsilon}(K(t))\le J_d^\varepsilon(I_\varepsilon)\le
\varepsilon^d\left(\frac{M(\varepsilon,t)}{\varepsilon}+N(\varepsilon,t)+1\right)\le
h(\varepsilon)\varepsilon^d\left(\frac{M(\varepsilon,t)}{h(\varepsilon)\varepsilon}+
\frac{N(\varepsilon,t)+1}{h(\varepsilon)}\right).$$ As $h$ is
regularly
 varying as $\varepsilon \rightarrow 0+$ with index $-\beta$, we
 deduce that for $d>\beta,
 h(\varepsilon)\varepsilon^d\rightarrow 0$ as $\varepsilon \rightarrow
 0+$ and so $H_d(K(t))=0$. This proves  that
 $\dim_{\mathcal{H}}K(t)\le \beta$.

\vspace*{0,5cm}

 Let us now prove the lower bound. We first prove the lower bound for a homogeneous $i$-fragmentation, i.e. we suppose here
 that $\alpha=0$.
 Let us fix $T_0>0$ and search for a lower bound of the Hausdorff dimension of
 $K(T_0)$. The two conditions of Hypothesis (H4) are symmetric by the transformation
$x\rightarrow 1-x$, so, without loss of generality, we suppose
here that $\liminf \frac{\nu(g\ge \varepsilon)}{ \nu(d\le
1-\varepsilon)}>0$. Hence there exists a constant $C$ such that for
$\varepsilon$ small enough we have $C\nu(g\ge \varepsilon)\ge
\nu(d\le 1-\varepsilon)$.
 We denote by $]g_t,d_t[$ the largest interval of the
 fragmentation at time $t$ and $T=\inf\{t\ge 0, d_t-g_t\le 1/2\}\wedge T_0$.
 So, for $0<s<t<T$, $]g_t,d_t[\subset]g_s,d_s[$. The idea is to
 prove that $\dim_{\mathcal{H}}\{g_t,0<t<T\}\ge \beta$ and as
 $\{g_t,0<t<T\}\subset K(T_0)$, we will conclude that
 lower bound holds for $\dim_{\mathcal{H}}  K(T_0)$.

 We know that $(U(t),t\ge 0)$ can be constructed from a PPP on
 $\mcr\times \ur \times \mcn$ with intensity measure $dt\times
 \nu\times \sharp$. So we have
 $$g_t=\sum_{s\in \mathcal{D}\cap[0,t]}\xi_s(d_{s^-} -g_{s^-}),$$
 where $(s,\xi_s)_{s\in \mathcal{D}}$ are the atoms of a Poisson
 measure on $\mcr \times [0,1]$ with intensity $ds\times \nu(g\in \cdot)$. We introduce now
 $$\sigma_t=\sum_{s\in \mathcal{D}\cap[0,t]}\xi_s.$$
 Then $\sigma$ is a subordinator with Levy measure $\Lambda(d\varepsilon)=\nu(g\in
 d\varepsilon)$ and we have :
 $$\forall \;0<s<t<T, \; g_t-g_s\ge \frac{1}{2}(\sigma_t-\sigma_s),$$
since $d_s-g_s< 1/2$ for $s\le T$.

It is then well known that, if we want to prove that $\dim_{\mathcal{H}}\{g_t,0<t<T\}\ge \gamma$, it is sufficient to prove that $g^{-1}$
 is Hölder-continuous
with exponent $\gamma$. We have then the following lemma :

\begin{duge}
Let $(f(t),0\le t\le T)$ and $(h(t),0\le t\le T)$ be two strictly  increasing càdlàg functions such that for all $0<s<t<T$, we have $h(t)-h(s)\ge
\frac{1}{2}(f(t)-f(s))$.
 Define
 $f^{-1}(x)=\inf\{u\ge 0, f(u)> x\}$ and suppose that $f^{-1}$ is Hölder-continuous with exponent $\gamma$.
Then $h^{-1}$ is also Hölder-continuous with exponent $\gamma$.
\end{duge}

Proof of the lemma. Let $s\ge t$ be  two elements of the set $H=\{h(t),0\le t\le T\}$. Hence there exist $x\ge y$ such that $h(x)=s$ and $h(y)=t$.
Then we have, for some constant $K$
$$h^{-1}(t)-h^{-1}(s)=y-x=f^{-1}\circ f(y)-f^{-1}\circ f(x)\le K(f(y)-f(x))^{\gamma}.$$
Besides we have $\quad t-s=h(y)-h(x)\ge \frac{1}{2}(f(y)-f(x)),$
so we get :
$$h^{-1}(t)-h^{-1}(s)\le 2^\gamma K (t-s)^{\gamma}.$$
Furthermore, $h^{-1}$ is constant on the interval components of $H^{c}$, and it follows then
$$h^{-1}(t)-h^{-1}(s)\le 2^\gamma K (t-s)^{\gamma} \quad \mbox{for all } s<t.\quad \Box$$

 Hence to prove that $\dim_{\mathcal{H}}\{g_t,0<t<T\}\ge \beta$, we just have to
 prove that $\sigma^{-1}$ is Hölder-continuous with exponent $\gamma$ for all $\gamma<\beta$.
 We use then the following lemma :

 \begin{duge}
\cite{Bertoin96} Let $(\sigma_s,s\ge 0)$ be a subordinator with no
drift and Lévy measure $\Lambda$. Let
$\Phi(\lambda)=\int_{0}^{\infty}(1-e^{-\lambda x})\Lambda(dx)$ and
$\gamma=\sup\{\alpha>0, \lim_{\lambda\rightarrow \infty}
\lambda^{-\alpha}\Phi(\lambda)=\infty\}$. Then, for every
$\varepsilon>0$, $\sigma^{-1}$ is a.s. Hölder-continuous
on compact intervals with exponent $\gamma-\varepsilon$.
 \end{duge}

To finish the proof of the homogeneous case, we have now to study
$\Lambda(d\varepsilon)=\nu(g\in d\varepsilon)$. In the following
we denote by $k$ an integer such that $\nu(s_k>0)=0$.

We remark that $\{g\ge \varepsilon\}\subset
\{Card\{i,s_i>\varepsilon/k\}\ge 2\}$, so $h(\varepsilon/k) \ge
\nu(g\ge \varepsilon)$. We notice also that $h(\varepsilon)\le
k\nu(g\ge \varepsilon \mbox{ or } d\le 1-\varepsilon)$. As
$\nu(d\le 1-\varepsilon)\le C\nu(g\ge \varepsilon)$ we get
$$\frac{h(\varepsilon)}{(C+1)k}\le \nu(g\ge \varepsilon)\le
h(\varepsilon/k).$$

Using that $h$  is regularly
 varying as $\varepsilon \rightarrow 0+$ with index $-\beta$, an
easy calculus proves that $\sup\{\alpha>0, \lim_{\lambda\rightarrow
\infty} \lambda^{-\alpha}\Phi(\lambda)=\infty\}=\beta$ and so
 $\sigma^{-1}$ is Holder-continuous with exponent $\beta-\varepsilon$ for all $\varepsilon>0$.
Hence we get that for each $t>0$, $\dim_{\mathcal{H}}K(t)=\beta$
a.s. As for $t<s$, $K(t)\subset K(s)$, $\dim_{\mathcal{H}}K(t)$
increases with $t$, and so we have also
$\dim_{\mathcal{H}}K(t)=\beta$ for all $t>0$ simultaneously a.s.

It remains now to prove the lower bound for an $i$-fragmentation
with strictly positive index  of self similarity . Let us use now
Proposition \ref{change} which changes the index
 of self-similarity of a fragmentation.
Let $(U^{\alpha}(t),t\ge 0)$ be a self similar fragmentation
fulfilling $(H)$. We write $U^{\alpha}(t)=U(T_t^{\alpha})$ as in
Proposition \ref{change} where $(U(t),t\ge 0)$ is a homogeneous
fragmentation. We denote by $(g_t,t\ge 0)$ (resp.
$(g^{\alpha}_t,t\ge 0)$) the left bound of the largest interval
component of $U(t)$ (resp. $U^{\alpha}(t)$). We know that for all
$T>0$, $\dim_{\mathcal{H}}\{g_t,0\le t\le T\}\ge \beta$. Or for $t$
small enough, we have $g_t^{\alpha}=g_{f(t)}$ where $f$ is a
continuous increasing function, so  for all $t>0$, there exists
$t'>0$ such that$$\dim_{\mathcal{H}}(K^{\alpha}(s),0\le s \le t)\ge
\dim_{\mathcal{H}}(g^{\alpha}_s,0\le s \le t)\ge
\dim_{\mathcal{H}}(g_s,0\le s \le t')\ge \beta.\; \Box$$

\begin{cor}\label{pack}
 Let $\nu$ be a dislocation measure fulfilling  Hypothesis \ref{hyp}.
Let $(U(t),t\ge 0)$ be a self-similar $i$-fragmentation with
characteristics $(\nu,0,0,\alpha)$   with $\alpha>0$. Let
$K(t)=[0,1]\backslash U(t)$. Then the packing dimension of $K(t)$ is
$\beta$
 for all $t>0$ simultaneously, a.s.
\end{cor}

Proof. Let us first recall the definition of the packing dimension
\cite{Tricot82}.  For  a subset $E\subset \mcr$ and $\alpha> 0$, let us define
$$M_{\alpha}(E)=\lim_{\varepsilon\rightarrow 0+}\sup\left\{\sum_{i=1}^{\infty} (2r_i)^\alpha, \,[x_i-r_i,x_i+r_i] \mbox{ disjoint}, x_i \in E,
 r_i<\varepsilon\right\},$$
and
$$\widehat{M}_\alpha(E)=\inf\left\{\sum_{n=1}^{\infty} M_\alpha(E_n),\; E\subseteq \bigcup_{n=1}^{\infty}E_n\right\}.$$
The packing dimension of $E$ is defined by
$$\dim_\wp(E)=\inf\{\alpha> 0, \; \widehat{M}_\alpha(E)=0\}=\sup\{\alpha> 0, \; \widehat{M}_\alpha(E)=\infty\}.$$

For  a subset $E\subset \mcr$   and $\varepsilon >0$, let $Z(E,\varepsilon)$ be the smallest number of interval of lengths
$2\varepsilon$ needed to cover $E$. We define
$$\Delta(E)=\limsup_{\varepsilon\rightarrow 0}\frac{\log Z(E,\varepsilon)}{-\log \varepsilon}.$$
Tricot \cite{Tricot82} proved that we have :
$$dim_\wp (E)=\inf\left\{\sup_n\Delta(E_n),\,E\subset \cup_n E_n\right\}.$$
It is then easy to see that for all $E\subset \mcr$,  we have $\dim_\mathcal{H} E\le \dim_\wp E$. Hence, to prove Corollary \ref{pack}, we just have to
get an upper bound of the packing dimension of $K(t)$. We use the same idea as for the Hausdorff dimension.
Let
$I_\varepsilon=]0,1[\backslash\{\mbox{ interval components of
}U(t)\mbox{ which size is larger than }\varepsilon\}$. So we have
$K(t)\subset I_{\varepsilon}$ and
$|I_{\varepsilon}|=M(\varepsilon,t)$ since $\nu$ is conservative.
Furthermore, $I_{\varepsilon}$ has at most $N(\varepsilon,t)+1$
interval components. We deduce that
$$Z(K(t),\varepsilon)\le Z(I_{\varepsilon},\varepsilon)\le h(\varepsilon)\left(\frac{M(\varepsilon,t)}{2\varepsilon h(\varepsilon)}
+\frac{N(\varepsilon,t)+1}{h(\varepsilon)}\right).$$
We get
$$\dim_\wp (K(t))\le\Delta(K(t))\le \limsup_{\varepsilon\rightarrow 0}\frac{\log h(\varepsilon)}{-\log \varepsilon}=\beta.$$
Hence, the packing dimension of the subset $K(t)$ coincides almost surely with its Hausdorff dimension (such subset is called "regular subset"). $\Box$

To conclude this section, let us discuss an example. We consider the
$m$-fragmentation introduced by Aldous and Pitman
\cite{Aldouspitman98} to study the standard additive coalescent.
Bertoin \cite{Bertoin00b} gave a construction of an $i$-fragmentation
$(U(t),t\ge 0)$ whose projection on $\srb$ is this
fragmentation. More precisely, let $\varepsilon=(\varepsilon_s,s\in
[0,1])$ be a standard positive Brownian excursion. For every $t\ge
0$, we consider
$$\varepsilon^{(t)}_s=ts-\varepsilon_s, \hspace*{1cm} S^{(t)}_s=\sup_{0\le u\le s
}\varepsilon^{(t)}_u.$$ We define $U(t)$ as the constancy intervals
of $(S^{(t)}_s,0\le s\le 1)$. Bertoin \cite{Bertoin02} proved also
that $(\freq{U(t)},t\ge 0)$ is an $m$-fragmentation with index of self
similarity $1/2$ and its dislocation measure  is carried by the
subset of sequences
$$\{s=(s_1,s_2,\ldots)\in \srb, s_1=1-s_2 \mbox{ and } s_i=0 \mbox{ for } i\ge
3\}$$ and  is given by  $$\tilde{\nu}_{AP}(s_1\in dx)=\left(2\pi x^3
(1-x)^3\right)^{-1/2}dx.$$ This proves that (H1), (H2) and (H3) hold with
$\beta=1/2$. Besides, as
$$\lim_{s\rightarrow 0^+}\frac{\varepsilon_s}{s}=\infty  \; a.s.,$$
0 is almost surely an isolated point of $[0,1]\backslash U(t)$  and
this implies that $\nu_{AP}(g>0)$ is finite. Hence we have $\limsup
\frac{\nu(g\ge \varepsilon)}{ \nu(d\le 1-\varepsilon)}<\infty$ and
Hypothesis (H3) holds. By Theorem \ref{theohaus}, we deduce that the
 Hausdorff  dimension of $[0,1]\backslash U(t)$ is $\frac{1}{2}$ a.s., a fact that can be checked directly using
 properties of Brownian motion.

\section{Interval components in uniform random order}

\begin{defi}\label{uniform}
Let $\tilde{\nu}$ be a measure on $\srb$ such that $\tilde{\nu}(\sum_i s_i<1)=0$. We define  $\widehat{\nu}$ as the measure
on $\ur$ which projection on $\srb$ is $\tilde{\nu}$ and which interval components are in uniform random order. More precisely,
set $s=(s_i)_{i\in\mcn}\in \srb$ with law $\tilde{\nu}$.
Let $(V_i)_{i\in\mcn}$ be  iid random variables uniform on $[0,1]$.
We denote then $U$ the random open subset of $]0,1[$ such that, if the
decomposition of U in disjoint open intervals ranked by their length
is $\coprod_{i=1}^\infty U_i$, we have
\begin{itemize}
\item For all $i\in \mcn$, $|U_i|=s_i$. \item For all $i\neq j,\;
U_i \prec U_j\;\Leftrightarrow \;V_i\le V_j$.
\end{itemize}
Since we have $\sum_i s_i =1$ a.s., there exists almost surely a
unique open subset of $]0,1[$ fulfilling this two points. We denote by
$\widehat{\nu}$ the distribution of $U$.
\end{defi}

\begin{prop}\label{unif}
Let  $(U(t),t\ge 0)$ is an $i$-fragmentation with measure $(\nu,0,0)$ and
such that for all $t\ge 0$, $U(t)$ has interval components in uniform random order. Then
$\nu$ has also interval components in uniform random order.
\end{prop}

Proof. Let $(F(t),t\ge 0)$ be the projection of $(U(t),t\ge 0)$ on $\srb$. We know that $F$ is then a $m$-fragmentation with measure
$(\tilde{\nu},0)$ where $\tilde{\nu}$ is the image of $\nu$ by the canonical projection $\ur\rightarrow \srb$.
Let $\gamma \in \sr_n$. Let $\pi\in \pn$ be the image of $\gamma$ by the canonical projection $\wp_1$ between $\sr$ and $\pinf$.
Let now remark that we have
$$q_\gamma=\frac{1}{s}\lim_{s\rightarrow 0}\mcp(\Gamma_{[n]}(s)=\gamma)=\frac{1}{k!}q_\pi,$$
where $k$ is the number of blocks of $\gamma$ and $q_\pi$ the jump rate of the $p$-fragmentation.
Let $\widehat{\nu}$ be the measure on $\ur$ obtained in Definition \ref{uniform} from $\nu$.
Let us recall that $\mathcal{Q}_{\infty,\gamma}=\{\gamma'\in \sr, \gamma'_{[n]}=\gamma\}$
and define also $\mathcal{P}_{\infty,\pi}=\{\pi'\in \pinf, \pi'_{[n]}=\pi\}$.
We have then
$$P^{\widehat\nu}(\mathcal{Q}_{\infty,\gamma})=\frac{1}{k!}P^{\tilde\nu}(\mathcal{P}_{\infty,\pi})=\frac{1}{k!}q_\pi=q_\gamma
=P^{\nu}(\mathcal{Q}_{\infty,\gamma}).$$
So we get that $\nu=\widehat{\nu}$ and hence $\nu$ has interval components in uniform random order. $\Box$

\vspace*{0.4cm}

Let us notice that the proof uses  $q_\gamma=\frac{1}{k!}q_\pi$, so if we want to extend this proposition to the time-inhomogeneous
case, we must  not only  suppose that $U(t)$ has interval components in uniform random order, but more generally that the semi-group of $U(t)$, $q_{t,s}(]0,1[)$ has
 interval components in uniform random order for all $t\ge 0$ and for all $s>t$.

\vspace*{0.4cm}

Conversely, we can ask if  $(U(t),t\ge 0)$ is an $i$-fragmentation with measure $(\nu,0,0)$ and
$\nu$ has  interval components in uniform random order, does this implies that $U(t)$ has interval components in uniform random order ?
The answer is clearly negative. Indeed, let $\nu$ be the following measure :
$$\nu=\delta_{U_1}+\delta_{U_2}
\mbox{ with } U_1=\Big]0,\frac{1}{3}\Big[\cup\Big]\frac{1}{3},\frac{2}{3}\Big[\cup\Big]\frac{2}{3},1\Big[
\mbox{ and  } U_2=\Big]0,\frac{1}{2}\Big[\cup\Big]\frac{1}{2},1\Big[.$$

Then $\nu$ has interval components in uniform random order, but
$U(t)$ has not this property since we have
$$\mcp\left(U(t)=\Big]0,\frac{1}{3}\Big[\cup\Big]\frac{1}{3},\frac{1}{2}\Big[\cup\Big]\frac{1}{2},\frac{2}{3}\Big[\cup\Big]\frac{2}{3},1\Big[\right)>0$$
and
$$\mcp\left(U(t)=\Big]0,\frac{1}{6}\Big[\cup\Big]\frac{1}{6},\frac{1}{2}\Big[\cup\Big]\frac{1}{2},\frac{5}{6}\Big[\cup\Big]\frac{5}{6},1\Big[\right)=0.$$

\subsection{Ruelle's fragmentation}\label{Ruelle}
In this section, we specify the semi-group of Ruelle's fragmentation
seen as an interval fragmentation. Let us recall the construction of
this interval fragmentation \cite{Bertoin00}.

 Let
$(\sigma^{*}_{t},0<t<1)$ be a family of stable subordinators such
for every $0<t_n<\ldots<t_1<1$,
$(\sigma^{*}_{t_1},\ldots,\sigma^{*}_{t_n})\overset{law}{=}
(\sigma_{t_1},\ldots,\sigma_{t_n})$ where
$\sigma_{t_i}=\tau_{\alpha_1}\circ\ldots \circ \tau_{\alpha_i}$ and
$(\tau_{\alpha_i},1\le i\le n)$ are $n$ independent stable
subordinators with indices $\alpha_1,\ldots,\alpha_n$ such that
$t_i=\alpha_1\ldots\alpha_i$. Fix $t_0\in ]0,1[$ and for $t\in
]t_0,1[$ define $T_{t}$ by :
$$ \sigma^{*}(T_t)=\sigma^*_{t_0}(1).$$
Then consider  the open subset :
$$U(t)=\Big]0,1\Big[\Big\backslash \left\{\frac{\sigma^{*}_{t}(u)}{\sigma^*_{t_0}(1)},0\le u\le T_{t}\right\}^{cl}.$$

 Bertoin and Pitman proved that
$(U(t),t\in[t_0,1[)$ is an $i$-fragmentation (with initial state
$U(t_0)\neq \un$ a.s.) and the semi-group of transition at time $t$
to time $s$ of the $m$-fragmentation $(\freq{U(t)},t\in [t_0,1[)$
 is $PD(s,-t)$-FRAG where $PD(s,-t)$ denotes the Poisson-Dirichlet law with parameter $(s,-t)$ (see
 \cite{Pitmanyor97} for more
 details about the Poisson-Dirichlet laws). Furthermore, the instantaneous dislocation measure of this $m$-fragmentation
 at time $t$ is $\frac{1}{t}PD(t,-t)$ (cf. \cite{Basdevant05}).
  We would like now to calculate the dislocation measure of
the $i$-fragmentation $(U(t),t\in[t_0,1[)$.

\begin{duge}
Let us define  $\widehat{PD}(t,0)$ as the measure on $\ur$ obtained
from $PD(t,0)$ by Definition \ref{uniform}. The distribution at time
$t$ of $U(t)$ is   $\widehat{PD}(t,0)$.
\end{duge}

Proof. For $t\in ]t_0,1[$, we have $\sigma^*_{t_0}=\sigma^*_t\circ
\tau_{\alpha}$ where $\alpha t=t_0$ and $\tau_{\alpha}$ is a stable
subordinator with index $\alpha$ and independent of $\sigma^*_t$.
Hence we get
$$U(t)=\Big]0,1\Big[\Big\backslash \left\{\frac{\sigma^{*}_{t}(u)}{\sigma^*_{t}(\tau_{\alpha(1)})},0\le u\le \tau_\alpha(1)\right\}^{cl}.$$
We can thus write
$$U(t)=\Big]0,1\Big[\Big\backslash\left\{\frac{\sigma_t(x)}{\sigma_t(a)},x\in[0,a[\right\}^{cl},$$
where $\sigma_t$ is a stable subordinator with index $t$ and $a$ is
a random variable independent of $\sigma_t$. If we denote by
$(t_i,s_i)_{i\ge 1}$ the time and size of the jump of $\sigma_t$ in
the interval $[0,a[$ ranked by decreasing order of the size of the
jumps, this family has the same law of $(t_{\tau(i)},s_{i})_{i\ge
1}$ for any $\tau$ permutation of $\mcn$. $\Box$

\begin{prop}
The semi-group of transition of the Ruelle's interval fragmentation
from time $t$ to time $s$ is $\widehat{PD}(s,-t)$-FRAG and the
instantaneous dislocation measure at time $t$  is
$\frac{1}{t}\widehat{PD}(t,-t)$.
\end{prop}

We would like now to apply Proposition  \ref{unif} to determine the
instantaneous measure of dislocation of Ruelle's fragmentation, but
this proposition holds only for time-homogeneous fragmentation. If
the fragmentation is inhomogeneous in time, we must first prove that
the semi-group of $U(t)$ has interval component in uniform order.
Fix $t\ge 0$ and $s>t$. Fix $y\in ]0,1[$ and denote by $I(t)$ the
interval component of $U(t)$ containing $y$. We shall prove that
$U(s)\cap I(t)$ has its interval component in uniform random order.
By the construction of $U(t)$, there exists $x\in]0,T_t[$ such that
$$I(t)=\Big]\frac{\sigma^{*}_{t}(x^-)}{\sigma^*_{t_0}(1)},\frac{\sigma^{*}_{t}(x)}{\sigma^*_{t_0}(1)}\Big[.$$
We have $\sigma^*_t=\sigma^*_s\circ \tau_{t/s}$ where $\tau_{t/s}$
is a stable subordinator with index $t/s$ and is independent of
$\sigma^*_{t+s}$. Hence, we get :
$$U(s)\cap I(t)=I(t)\Big\backslash
\left\{\frac{\sigma^{*}_{s}(y)}{\sigma^*_{t_0}(1)},\;\tau_{t/s}(x^-)\le
y\le \tau_{t/s}(x)\right\}^{cl}.$$ Since $\tau_{t/s}$ is independent
of $\sigma^*_s$, the jump of $\sigma^*_s$ on the interval
$]\tau_{t/s}(x^-), \tau_{t/s}(x)[$ are in uniform random order.
Since as $m$-fragmentation the semi-group of transition is
$PD(s,-t)$-FRAG, we deduce that, as $i$-fragmentation, the
semi-group is $\widehat{PD}(s,-t)$-FRAG. To prove that the
dislocation measure at time $t$ is $\frac{1}{t}\widehat{PD}(t,-t)$,
we just have to apply the Proposition \ref{unif}. $\Box$

\subsection{Dislocation measure of the fragmentation derived from the additive coalescent}
Recall the construction of an $i$-fragmentation $(U(t),t\ge 0)$
 from a Brownian motion exposed in Section \ref{Haussec}. We already know its characteristics as a
$m$-fragmentation : the erosion rate is null, the index of self similarity is equal to $1/2$ and the dislocation measure $\tilde{\nu}_{AP}$ is given
by :
$$ \tilde{\nu}_{AP}(s_1\in dx)=(2\pi x^3(1-x^3))^{-1/2}dx \quad \mbox{ for } x\ge 1/2, \qquad   \tilde{\nu}_{AP}(s_1=1-s_2)=1.$$

\begin{prop}
The $i$-fragmentation derived from a Brownian motion \cite{Bertoin00b} has  dislocation measure
$\nu_{AP}$ such that :
\begin{itemize}
\item $\nu_{AP}$ is carried by the subset of ]0,1[ shaped as $]0,1[ \setminus x$. So we will write $\nu_{AP}(x)$ instead of $\nu_{AP}(]0,1[\setminus x)$.
\item For all $x\in ]0,1[, \quad \nu_{AP}(dx)=(2\pi x(1-x^3))^{-1/2}dx$.
\end{itemize}

\end{prop}

Notice that we have $\nu_{AP}(dx)=x\tilde{\nu}_{AP}(s_1\in dx \mbox{
or }s_2\in dx)$ for all $x \in ]0,1[$. Hence,  given that the
$m$-fragmentation splits in two block of size $x$ and $1-x$, the
left block of the $i$-fragmentation will be a size biased pick of
$x$ and $1-x$.

Proof. The first part of the proposition is immediate since we have $ \tilde{\nu}_{AP}(s_1=1-s_2)=1$. For the second part,
let us use Theorem 9 in \cite{Bertoin00b} which gives the distribution $\rho_t$ of the most left fragment of $U(t)$ :
$$\rho_t(dx)=t\frac{1}{\sqrt{2\pi x(1-x)^3}}\exp\left(-\frac{xt^2}{2(1-x)}\right)dx\quad \mbox{ for  all } x\in ]0,1[.$$

We get
$$\nu_{AP}(dx)=\lim_{t\rightarrow 0}\frac{1}{t}\rho_t(dx)=\frac{1}{\sqrt{2\pi x(1-x)^3}}. \,\Box$$

\vspace*{0.5cm}

We can also give a description of the distribution at time $t>0$ of
$U(t)$. Recall the result obtained by Chassaing and Janson
\cite{Chassaingjanson01}. For a random process $X$ on $\mcr$ and
$t\ge 0$, we define $\ell_t(X)$ as the local time of $X$ at level 0
on the interval $[0,t]$, i.e.
$$\ell_t(X)=\lim_{\varepsilon\rightarrow
0+}\frac{1}{2\varepsilon}\int_{0}^{t}\un_{\{|X_s|<\varepsilon\}}ds,$$
whenever the limit makes sense.

Let $X^t$ be a reflected Brownian bridge conditioned on
$\ell_1(X^t)=t$. We define $\mu\in]0,1[$ such that
$$\ell_\mu(X^t)-t\mu=\max_{0\le u\le 1}\ell_u(X^t)-tu.$$
It is well known that this equation has almost surely a unique
solution. Let us define the process $(Z^t(s),0\le s\le 1)$ by
$$Z^t(s)=X^t(s+\mu \,[\mbox{mod }1]).$$
Chassaing and Janson \cite{Chassaingjanson01} have proved that for
each $t\ge0$
$$U(t)\overset{law}= ]0,1[\backslash\{x\in [0,1], Z^t(x)=0\}.$$
Besides, as the inverse of the local time of $X^t$ defined by
$$T_x=\inf\{u\ge 0,\, \ell_u(X^t)>x\}$$
is a stable subordinator with Lévy measure $(2\pi x^3)^{-1/2}dx$
conditioned to $T_t=1$, we deduce the following description of the
distribution of $U(t)$ :

\begin{cor}Let $t>0$.
Let $T$ be a stable subordinator with Lévy measure $(2\pi
x^3)^{-1/2}dx$ conditioned to $T_t=1$. Let us define $m$ as the
unique number  on $[0,t]$ such that
$$tT_{m^-}-m\le tT_u-u\quad \mbox{ for all } u\in [0,t],$$
where $T_{m^-}=\lim_{x\rightarrow m^-}T_x$. We set :
$$\begin{array}{ccll}
\tilde{T}_x
&=&T_{m+x}-T_{m^{-}}& \mbox{ for } \,0< x < t-m, \\
&& T_{m+x-t}-T_{m^{-}}+1 &\mbox{ for }\, t-m\le x \le t.
\end{array}$$
Then $$U(t)\overset{law}=]0,1 [\backslash\{\tilde{T}_x,x\in
[0,t]\}^{cl}.$$
\end{cor}

Proof. It is clear that $\{u, X^t(u)=0\}$ coincides with $\{T_x,x\in
[0,t]\}^{cl}$ when $T$ is the inverse of the local time of $X^t$.
Hence, we just have to check that if we set $m=\ell_{\mu}(X^t)$,
then $m$ verifies the equation $tT_{m^-}-m\le tT_u-u$  for all $u\in
[0,t]$. Since $X^{t}(\mu)=0$, we have $T_{m^-}=\mu$, thus we get :
$$tT_{m^-}-m=t\mu-\ell_{\mu}(X^t)\le tv-\ell_v(X^t) \mbox{ for all }
v\in [0,1].$$ Let us fix $u\in [0,t]$. Since $\ell_v(X^t)$ is a
continuous function, there exists $v\in [0,1]$ such that
$\ell_v(X^t)=u$. Besides we have $T_u^{-}\le v\le T_u$, so we get
 $$tT_{m^-}-m\le tT_u-u. \; \Box$$

Hence, the distribution of $[0,1]\setminus U(t)$ can be obtained  as
the closure of the range of a stable subordinator $(T_s,0\le s \le
t)$ with index $1/2$ and conditioned on $T_t=1$ randomly shifted
(recall also that Chassaing and Jason \cite{Chassaingjanson01} have
proved that the left most fragment of $U(t)$ is  size-biased
picked).

 \vspace*{1cm} \nocite{Haasmiermont04}

\bibliographystyle{plain}
\bibliography{biblio}

\begin{thebibliography}{10}

\bibitem{Aldouspitman98}
D.~Aldous and J.~Pitman.
\newblock The standard additive coalescent.
\newblock {\em Ann. Probab.}, 26(4):1703--1726, 1998.

\bibitem{Basdevant05}
A.-L. Basdevant.
\newblock Ruelle's probability cascades seen as a fragmentation process.
\newblock Preprint, 2005.

\bibitem{Berestycki02}
J.~Berestycki.
\newblock Ranked fragmentations.
\newblock {\em ESAIM}, 6:157, 2002.

\bibitem{CoursBertoin03}
J.~Bertoin.
\newblock {\em Random fragmentation and coagulation processes}.
\newblock Cambridge University Press, Cambridge.
\newblock To appear.

\bibitem{Bertoin96}
J.~Bertoin.
\newblock {\em Lévy processes}.
\newblock Cambridge University Press, Cambridge, 1996.

\bibitem{Bertoin00b}
J.~Bertoin.
\newblock A fragmentation process connected to brownian motion.
\newblock {\em Probab. Theory Related Fields}, 117(2):289--301, 2000.

\bibitem{Bertoin02}
J.~Bertoin.
\newblock Self-similar fragmentations.
\newblock {\em Ann. Inst. H. Poincar\'e Probab. Statist.}, 38(3):319--340,
  2002.

\bibitem{Bertoin04}
J.~Bertoin.
\newblock On small masses in self-similar fragmentations.
\newblock {\em Stochastic Process. Appl.}, 109(1):13--22, 2004.

\bibitem{Bertoinlegall00}
J.~Bertoin and J.-F. Le~Gall.
\newblock The {B}olthausen-{S}znitman coalescent and the genealogy of
  continuous-state branching processes.
\newblock {\em Probab. Theory Related Fields}, 117(2):249--266, 2000.

\bibitem{Bertoin00}
J.~Bertoin and J.~Pitman.
\newblock Two coalescents derived from the ranges of stable subordinators.
\newblock {\em Electron. J. Probab.}, 5:no. 7, 2000.

\bibitem{Bolthausensznitman98}
E.~Bolthausen and A.-S. Sznitman.
\newblock On \uppercase{R}uelle's probability cascades and an abstract cavity
  method.
\newblock {\em Comm. Math. Phys.}, 197(2):247--276, 1998.

\bibitem{Chassaingjanson01}
P.~Chassaing and S.~Janson.
\newblock A {V}ervaat-like path transformation for the reflected {B}rownian
  bridge conditioned on its local time at 0.
\newblock {\em Ann. Probab.}, 29(4):1755--1779, 2001.

\bibitem{Falconer86}
K.~J. Falconer.
\newblock {\em The geometry of fractal sets}, volume~85 of {\em Cambridge
  Tracts in Mathematics}.
\newblock Cambridge University Press, Cambridge, 1986.

\bibitem{Gnedin97}
A.~V. Gnedin.
\newblock The representation of composition structures.
\newblock {\em Ann. Probab.}, 25(3):1437--1450, 1997.

\bibitem{Haasmiermont04}
B.~Haas and G.~Miermont.
\newblock The genealogy of self-similar fragmentations with negative index as a
  continuum random tree.
\newblock {\em Electron. J. Probab.}, 9:no. 4, 57--97 (electronic), 2004.

\bibitem{Kingman82}
J.F.C. Kingman.
\newblock The coalescent.
\newblock {\em Stochastic Process. Appl.}, 13:235--248, 1982.

\bibitem{Pitmanyor97}
J.~Pitman and M.~Yor.
\newblock The two-parameter {P}oisson-{D}irichlet distribution derived from a
  stable subordinator.
\newblock {\em Ann. Probab.}, 25:855--900, 1997.

\bibitem{Ruelle87}
D.~Ruelle.
\newblock A mathematical reformulation of \uppercase{D}errida's \uppercase{REM}
  and \uppercase{GREM}.
\newblock {\em Commun. Math. Phys.}, 108:225--239, 1987.

\bibitem{Tricot82}
C.~Tricot.
\newblock Two definitions of fractional dimension.
\newblock {\em Math. Proc. Cambridge Philos. Soc.}, 91(1):57--74, 1982.

\end{thebibliography}

\end{document}